\documentclass[a4paper,english,11pt]{article}

\usepackage[english]{babel}
\usepackage{theorem,epsfig,subfigure,array,fancyheadings}
\usepackage{varioref,float,amssymb,amsmath,float}
\usepackage[latin1]{inputenc}
\usepackage[active]{srcltx}

\floatstyle{ruled}
\newfloat{algorithm}{H}{alg}
\floatname{algorithm}{Algorithmus}

\newtheorem{satz}{Theorem}
\newtheorem{lemma}[satz]{Lemma}
\newtheorem{folge}[satz]{Corollary}

\newtheorem{rem}[satz]{Remark}

\newcommand{\N}{\mathop{\it I\!\!\:\! N}\nolimits}       
\newcommand{\R}{\mathop{\it I\!\!\!\, R}\nolimits}       
\newcommand{\eor}{\hspace*{\fill} $\square $}            
\newcommand{\eop}{\hspace*{\fill} $\blacksquare$}        
\newcommand{\bew}{\textbf{Proof:}}

\newcommand{\eps}{\varepsilon}

\newcommand{\nach}{\; \longrightarrow \;}
\newcommand{\auf}{\; \longmapsto \;}

\newcommand{\D}{\displaystyle}
\newcommand{\T}{\textstyle}
\newcommand{\la}{\lambda}

\newcommand{\ttilde}[1]{\raisebox{.15ex}{$\stackrel{\mbox{}_\approx}{{#1}}$}}

\renewcommand{\thefootnote}{\fnsymbol{footnote}} 
\renewcommand{\baselinestretch}{1}

\begin{document}
\thispagestyle{empty}
\mbox{} \vskip3cm
\centerline{\LARGE\bf Optimal exploitation of renewable resource stocks:}
\centerline{\LARGE\bf Necessary conditions}
\vskip1cm

\setcounter{footnote}{1}
\noindent
J. Baumeister\footnote{Corresponding author} \ Fachbereich Mathematik,
Johann Wolfgang Goethe -- Universit\"at, Robert Mayer Strasse 6 -- 10,
D 60054 Frankfurt/Main, Germany \\
{\tt e-mail:~Baumeister@math.uni-frankfurt.de} \\[1ex]
A. Leit\~ao \ Radon Institute for Computational and Applied Mathematics
Austrian Academy of Sciences, Altenbergerstr. 69, A-4040 Linz, Austria \\
{\tt e-mail:~antonio.leitao@oeaw.ac.at}
\vskip1cm

\begin{center}
\today 
\end{center}

\newpage
\thispagestyle{empty}
\mbox{} \vskip3cm
\centerline{\LARGE\bf Optimal exploitation of renewable resource stocks:}
\centerline{\LARGE\bf Necessary conditions}
\vskip1cm

\centerline{\Large\bf Summary}
\medskip
\noindent
We study a model for the exploitation of renewable stocks developed in
\cite{CC}. In this particular control problem, the control law contains a
measurable and an impulsive control component. We formulate Pontryagin's
maximum principle for this kind of control problems, proving first order
necessary conditions of optimality.
Manipulating the correspondent Lagrange multipliers we are able to define
two special switch functions, that allow us to describe the optimal
trajectories and control policies nearly completely for all possible
initial conditions in the phase plane.
\bigskip\bigskip


\newpage
%
\renewcommand{\thefootnote}{\arabic{footnote}}
\setcounter{footnote}{0}
\pagestyle{headings}
\pagenumbering{arabic}
\setcounter{page}{3}
%
%
%
\section{Introduction} \label{sec:intro}

\subsection{Description of the model}

Consider a bio-economic model of the commercial fishery under sole
ownership. The model is governed by the quantities described in
Table~\ref{tab:main-var}. \medskip

\begin{table}[h]
\begin{center}
\begin{tabular}{l@{\ \ :\ }l}
 \hline
 $t$    & time variable \\
 $x(t)$ & population biomass at time $t$ \\
 $h(t)$ & harvest rate at time $t$ \\
 $E(t)$ & fishing effort at time $t$ \\
 $K(t)$ & amount of capital invested in the fishery at time $t$ \\
 $I(t)$ & investment rate at time $t$ \\
 \hline
\end{tabular}
\end{center}
\vskip-0.4cm
\caption{Main relevant variables.} \label{tab:main-var}
\end{table}
The model is based on the following assumptions:
\begin{itemize}
\item $h(t) = q E(t) x(t)$; $q$ is a catch coefficient;
\item $x'(t) = F(x(t)) - q E(t) x(t)$, $t \ge 0$, $x(0) = x_0$;
      $F$ is the natural growth function;
\item $K'(t) = -\gamma K(t) + I(t)$, $t \ge 0$, $K(0) = K_0$;
      $\gamma \ge 0$ is the rate of depreciation;
\item constraints:
      $0 \le x(t),\ K(t),\ E(t)$; $E(t) \le K(t)$;
\item non-malleability: $0 \le I(t) \le \infty$, $t \ge 0$;
\item existence of two biological equilibrium points:
      $F(0) = F(\bar{x}) = 0$, $\bar{x} > 0$;
\item properties of the production function:
      $$ F \in C^2[0,\infty)\, ,\ F(x)>0\, ,\ 0 < x < \bar{x}\, ,\
         F''(x) < 0\, ,\ 0 \le x \le \bar{x}\, ; $$
\item objective function (discounted cash flow):
      $\int\limits_0^{\infty}e^{-\delta t}\{ph(t)-cE(t)-rI(t)\} dt$; \\
      $\delta > 0$ is the instantaneous rate of discount, $p \ge 0$ is the
      price of landed fish, $c \ge 0$ is the operating cost per unit effort,
      $r \ge 0$ is the price of capital.
\end{itemize}

A concrete production function satisfying the assumption above is given
by the logistic mapping $F(x) := ax (1-\frac{\D x}{\D k})$ (with $a>0$,
$k>0$). In our figures we use this production function.
%
%
%
\subsection{The optimal control problem}

We set $E = u \, K$ and consider $u$ as a second control variable. Without
loss of generality we use $q=1$. After this manipulation, the problem we
want to consider is the following optimal control problem: \\[2ex]
\noindent $Q(x_0,K_0) \hskip1cm
\left\{ \begin{array}{l}
\text{Minimize } J(x_0,K_0;I,u) :=\int\limits_0^{\infty}e^{-\delta t}
                 \{ rI(t)+cu(t)K(t)-pu(t)K(t)x(t) \} dt \\
\text{subject to} \\
\hspace{3ex} x'=F(x)-u(t) K(t)x\;,t \ge 0,\;x(0)=x_0\,,\\
\hspace{3ex} K'=-\gamma K+I(t)\;, t \ge 0,\,K(0) = K_0\,,\\
\hspace{3ex} 0 \le x(t), K(t)\, ,\ 0 \le u(t) \le 1\, ,\ I(t)\ge 0\, ,
             \ t\in[0,\infty)\, . \end{array} \right. $ \medskip

This problem is considered in \cite{CC} along the \glqq royal road\grqq of
Carath\'eodory. But the analysis is not rigorous in the details (see Section~%
\ref{sec_concur}). In \cite{Mu} the problem is given as an illustration
for the problem of the type considered in the paper but the results are not
applicable (for proving existence).

In \cite{VP} there is also a hint to this problem and finally, 
we find the problem in \cite{C2} considered as an example for the application 
of the maximum principle, but nothing has been made rigorous. For the 
background in fishery management see \cite{Fl}, \cite{An} and \cite{SS}.

It is known, see e.g. \cite{Mu}, that control problems may have no solution
if the control variable is unbounded and both the cost functional and the
dynamics depend linearly on the control. This situation is given here with
respect to the control $I$. In order to avoid non--existence we have to
replace the conventional control $I$ by an impulse control, i.e. jumps in
the state are allowed. Therefore, we consider $I$ as a Borel measure and
the capital function $K$ as a function of bounded variation. Then, problem
$Q(x_0,K_0)$ becomes:
\medskip

\noindent $P(x_0, K_0) \hskip1cm \left\{ \begin{array}{l}
\text{Minimize } J(x_0,K_0;\mu,u) := \int\limits_0^{\infty}e^{-\delta t} r
                \mu(dt) + \int_0^{\infty} e^{-\delta t} \{c - px(t)\}
                u(t)K(t) dt \\
\text{ subject to }
\hspace{3ex}        (u,\mu) \in U_{ad} \times C^*\, \text{ and } \\
\hspace{3ex}        x' = F(x)-u(t) K x\,,\ x(0)=x_0\, , \\
\hspace{3ex}        dK = -\gamma K dt + \mu(dt)\, ,\ K(0) = K_0\, . 
\end{array} \right. $ \medskip

\noindent  Here
\begin{eqnarray*}
U_{ad} &:=& \{v \in L_{\infty}[0,\infty)| 0 \le v(t) \le 1 \text{ a.e. in }
            [0,\infty)\} \, , \\
C^*    &:=& \{ \mu | \mu \text{ a non--negative Borel measure on }
            [0,\infty)\} \, .
\end{eqnarray*}
Notice that the constraints $0 \le x(t), K(t)$, $t \in [0,\infty)$, are
satisfied due to the assumptions above if $x_0 \ge 0$, $K_0 \ge 0$.

The initial value problem
\begin{equation}\label{eq:kdiff}
dK = -\gamma K dt + \mu(dt)\, ,\ K(0) = K_0
\end{equation}
has to be considered as differential equation with a measure: a function
$K : [0,t_1) \nach \R$ ($t_1 \in (0,\infty]$) is a solution if
\begin{equation} \label{eq:kint}
K(t) = K_0 - \int_0^t \gamma K(s) ds + \int_{[0,t]}\mu(ds)\, ,
       \ 0 \le t < t_1\, ;
\end{equation}
This implies that $K$ is right continuous in $(0,t_1)$ and $K(0) = K_{0,+}
\mu(\{t_0\})$ where $K_{0,+}$ denotes $\lim_{t \downarrow 0} K(t)$.

\begin{rem} \label{rem:sol}
Due to the fact that the coefficients in front of the control measure
$\mu$ does not depend on the state we may use the solution concept as
given above, the so called Young solutions (see \cite{Ri}). Otherwise
we would have to use the concept of robust solutions considered in
\cite{SV}, \cite{DR}, \cite{MR}, \cite{SR}, \cite{BL}.
\eor \end{rem}

We set $\kappa := \delta + \gamma$, $r' := r \kappa$, $c_* := c + r'$ and
define functions $g,\psi,\psi_*$ on $(0,\bar{x})$ by
\begin{eqnarray*}
g(x)      &:=& \delta - F'(x) + \frac{\D F(x)}{\D x}\, , \\
\psi(x)   &:=& (px-c)(\delta - F'(x))-\frac{\D c F(x)}{\D x}\, , \\
\psi_*(x) &:=& (px-c_*)(\delta - F'(x))-\frac{\D c_* F(x)}{\D x}\, .
\end{eqnarray*}
Further we consider the following conditions:
\begin{itemize}
\item [($V1$)~] $F \in C^2[0,\infty) \cap C^3(0,\bar{x})$; \\[1ex]
                $F(0) = F(\bar{x}) = 0\,;\,F(x)>0\,,\,0 < x < \bar{x}\, ;
                \, F''(x)<0\,,\,0 \le x \le \bar{x}$.
\item [($V2$)~] $\delta > 0$, $r > 0$, $c > 0$, $\gamma > 0$.
\item [($V3$)~] $c_* - p \bar{x} < 0$.
\item [($V4$)~] There exist $\tilde{x}$, $x^* \in (0,\bar{x})$ with
         $$\psi(x)<0\, ,\ 0 < x <\tilde{x}\, ,\ \psi(\tilde{x})= 0\, ,
                 \ \psi(x)>0\, ,\ \tilde{x} < x < \bar{x}\, , $$
         $$\psi_*(x)<0\, ,\ 0 < x < x^*\, ,\ \psi_*(x^*) = 0\, ,
                 \ \psi_*(x)>0\, ,\ x^* < x < \bar{x}\, . $$
\item [($V5$)~] $\psi'(x) > 0$, $x \in (0,\tilde{x})$, $\psi_*'(x) > 0$,
                $x \in (\tilde{x},\bar{x})$.
\item [($V6$)~] $g'(x) > 0$, $x \in (0,\bar{x})$.
\end{itemize}

\begin{rem} \label{rem:red}
Notice that the conditions $(V1)$, \dots , $(V6)$ are satisfied for the
logistic production function if the constants are chosen appropriately.
Notice too that $(V4)$, $(V5)$ contain redundant information.
\eor \end{rem}

\begin{rem}\label{rem:red1}
In the subsequent analysis it is very important that $x = \bar{x}$ is an
attracting equilibrium point, while $x = 0$ is an unstable equilibrium point.
\eor \end{rem}

At this point we define $\widetilde{K} := F(\tilde{x}) / \tilde{x}$\, and\,
$K^* := F(x^*) / x^*$. Due to assumption (V1), follows $g(x) > 0$ in
$[0,\bar{x}]$. Therefore, we have $\psi(x) > \psi_*(x)$ and, consequently,
$\tilde{x} < x^*$, $\widetilde{K} > K^*$ must hold.
\medskip

Under the assumptions $(V1)$, \dots , $(V6)$ one can prove existence of
optimal solutions of $P(x_0,K_0)$; see, e.g., \cite{Si}
%
%
%
\subsection{The verification approach} \label{sec_concur}

As already mentioned the problem $P(x_0,K_0)$ is considered in
\cite{CC}. By using a Hamilton--Jacobi--Bellman equation on $(0,\bar{x})
\times (0,\infty)$ a candidate for an optimal control policy is defined for
each $(x_0,K_0) \in (0,\bar{x}) \times (0,\infty)\,$. This results in the
definition of a function
$$ S : (0,\bar{x}) \times (0,\infty) \nach \R \, , $$
such that for all $(x,K) \in  (0,\bar{x}) \times (0,\infty)$, for all
$u \in [0,1]$ and for all $I \ge 0$ -- controls with jumps are avoided
by considering them as \glqq limits\grqq of regular controls --
\begin{multline} \label{eq:hjbccm}
\delta S(x,K)+F(x)S_x(x,K)-\gamma KS_K(x,K) \\
\ge I(S_K(x,K)+r) + uK \{qxS_x(x,K) - pqx + c\}
\end{multline}
holds. Then it is stated that $S$ is the value function $V$ where the
value function $V$ is given here by
$$ V(x_0,K_0) := \inf \{J(x_0,K_0;I,u)|(I,u) \text{ admissible}\}\, ,\
   (x_0,K_0) \in (0,\bar{x}) \times (0,\infty) \, . $$
Implicitly there are only used controls which result in states $x$, $K$
such that
$$ \int_{t \to \infty} e^{-\delta t} S(x(t),K(t)) = 0 \, . $$
One can follow the analysis in \cite{CC} partly but for some steps the
assumptions are not sufficient and some arguments are not complete. Since
$S$ is not differentiable everywhere they use the argument that each
problem $P(x_0, K_0)$ may be approximated by the problem $Q(x_0, K_0)$.
This density argument is a very deep topological argument and no results
to make this argument rigorous are available from the literature. It is
open whether on this road the verification of optimal controls is possible.
Thus, the verification of the optimal policy in \cite{CC} has to be
considered as an open problem.Two different ways may be considered in
order to circumvent these difficulties: Firstly, extension of the
so--called Hamilton--Jacobi--Bellman equation such that jumps are
allowed. Secondly, proof of the closure property inherent in the
density argument.
%
%
%
\subsection{Outline of the paper} \label{ssec:outline}

In Section~\ref{sec:nec} we study the necessary conditions, furnished by a
special version of Pontryagin's maximum principle (see Appendix). Through
manipulation of the Lagrange multipliers we manage to define two special
switch functions. The first switch helps to determine the bang--bang
behavior of the measurable component of the control policy, while the
second one is a jump switch, which gives us a necessary condition for
discontinuities in the state variables.

In Section~\ref{sec:exploit} we use the necessary conditions of optimality,
rewritten for the switch functions, in order to detect both optimal and non
optimal behavior of the admissible processes. Following trajectories
backwards in time and observing the evolution of the switches, we are able
to construct auxiliary curves in the phase plane $(x,K)$, that are very
useful to determine optimal behavior. Particularly we are able to detect
two jump curves in the phase plane. This shows that the application of the
Pontryagin maximum principle can be used in order to construct the extremals
of the problem.

In Section~\ref{sec:optim} we put all arguments together and summarize the
obtained results in the form of Theorem~\ref{th:rrr}. The next two theorems,
\ref{th:gamma} and \ref{th:sigma0}, treat some special initial conditions,
which may occur. However, the argumentation follow the spirit of
Theorem~\ref{th:rrr}.

It is worth to mention that our results are in agreement with the
conclusions in \cite{CC}.
%
%
%
\subsection{Interpretation of the main results} \label{ssec:mres}

In this section we provide the economic interpretation of our main result,
which is obtained in Theorem~\ref{th:rrr} and auxiliary Theorems~%
\ref{th:gamma} and~\ref{th:sigma0}. These theorems describe optimal
behavior of processes and corresponding Lagrange multipliers for all
initial conditions in the state space $[0,\bar{x}] \times [0,\infty)$.

For details on the notation, particularly the definition of the curves
$\Sigma^*$, $\widetilde{\Sigma}$, $\Sigma_s$, $\Sigma_0$, $\Gamma_1$,
$\Gamma_2$, $\Gamma_3$, $\Gamma_4$, the reader should refer to
Section~\ref{sec:exploit} (see also to Figure~\ref{fig:szenario}). The
main regions (R1), \dots, (R5) are defined in Section~\ref{sec:optim}
and are illustrated in Figure~\ref{fig:regions}.
\medskip

\noindent Case~1: $(x_0,K_0) = (x^*,K^*)$ \vspace{-0.2cm}
\begin{itemize} \item[]
One should invest with constant rate ($\mu = \gamma K^* dt$) and fish with
maximal effort ($u = 1$) for all $t \ge 0$. Consequently, the optimal
trajectory satisfies $(x(t),K(t)) = (x^*,K^*)$ for all $t \ge 0$ (this is
the first singular arc).
\end{itemize}

\noindent Case~2: $(x_0,K_0) \in \Sigma^*$ \vspace{-0.2cm}
\begin{itemize} \item[]
At the initial time $t=0$ one should make an impulsive investment, in
such a way that $(x_0, K_{0,+}) = (x^*,K^*)$. Then one should proceed
as in Case~1.
\end{itemize}

\noindent Case~3: $(x_0,K_0) \in$ (R2) \vspace{-0.2cm}
\begin{itemize} \item[]
One should not invest ($\mu = 0$) and fish with maximal effort ($u = 1$)
until the optimal trajectory reaches the curve $\Sigma^*$. Then one should
proceed as in Case~2.
\end{itemize}

\noindent Case~4: $(x_0,K_0) \in$ (R1) \vspace{-0.2cm}
\begin{itemize} \item[]
At the initial time $t=0$ one should make an impulsive investment, in
such a way that $(x_0, K_{0,+}) \in \Sigma_s$, the so called {\em jump
curve}. Then one should proceed as in Case~3.
\end{itemize}

\noindent Case~5: $(x_0,K_0) \in$ (R3) \vspace{-0.2cm}
\begin{itemize} \item[]
One should not invest ($\mu = 0$) and should not fish ($u = 0$) until the
optimal trajectory reaches the curve $\Sigma_0$ (this configures (R3) as
a {\em moratorium region}). Then one should proceed as in Case~3.
\end{itemize}

\noindent Case~6: $(x_0,K_0) \in \widetilde{\Sigma}$ \vspace{-0.2cm}
\begin{itemize} \item[]
One should not invest ($\mu = 0$) and should fish with moderate effort
$u(t) = K(t)^{-1} F(\tilde{x}) \tilde{x}^{-1}$ until the optimal trajectory
reaches the state $(\tilde{x}, \ttilde{K})$ (second singular arc). Note
that the fish population remains constant ($x = \tilde{x}$) during this
first time interval. Then one should proceed as in Case~3.
\end{itemize}

\noindent Case~7: $(x_0,K_0) \in$ (R4) \vspace{-0.2cm}
\begin{itemize} \item[]
One should not invest ($\mu = 0$) and should not fish ($u = 0$) until the
optimal trajectory reaches the curve $\widetilde{\Sigma}$ (this configures
(R4) as a {\em moratorium region}). Then one should proceed as in Case~6.
\end{itemize}

\noindent Case~8: $(x_0,K_0) \in$ (R5) \vspace{-0.2cm}
\begin{itemize} \item[]
One should not invest ($\mu = 0$) and should fish with maximal effort
($u = 1$) until the optimal trajectory reaches the curve $\widetilde{\Sigma}$.
Then one should proceed as in Case~6.
\end{itemize}

\noindent Case~9: $(x_0,K_0) \in \Sigma_0 \cup \Sigma_s \cup \Gamma_3$
\vspace{-0.2cm}
\begin{itemize} \item[]
This case is analog to Case~3.
\end{itemize}

\noindent Case~10: $(x_0,K_0) \in \Gamma_4$ \vspace{-0.2cm}
\begin{itemize} \item[]
This case is analog to Case~7.
\end{itemize}
%
%
%
\section{Necessary conditions}\label{sec:nec}

In this section we use the maximum principle to derive first order necessary
conditions for problem $P(x_0, K_0)$ and define, with the aid of the Lagrange
multipliers, two auxiliary functions (switches) that play a key rule in the
analysis of the optimal trajectories. We start defining the Hamilton function
$\widetilde{H}$ by
$$  \widetilde{H}(t, \check{x}, \check{K}, w, \widetilde{\la}_1,
    \widetilde{\la}_2,\eta) :=
    \widetilde{\la}_1 (F(\check{x}) - w \check{K} \check{x}) -
    \widetilde{\la}_2 \gamma \check{K} - \eta e^{-\delta t}
    (c - p \check{x}) w \check{K} \, . $$
Let $(u,\mu)$ be an optimal control policy for $P(x_0, K_0)$ and let $(x,K)$
be the associated state. From the maximum principle in the appendix we obtain
constants $\widetilde{\la}_{1,0}$, $\widetilde{\la}_{2,0}$, $\eta \in \R$ and
adjoint functions $\widetilde{\la}_1$, $\widetilde{\la}_2: [0,\infty) \to \R$
such that
\begin{eqnarray*}
\widetilde{\la}_{1,0}^2 + \widetilde{\la}_{2,0}^2 + \eta^2 &\neq&
       0\, ,\ \eta \ge 0\, ,\\
x' &=& F(x)-u(t) K x\, ,\ x(0) = x_0\, , \\
dK &=& -\gamma K dt + \mu(dt)\, ,\ K(0) = K_0\, , \\
\widetilde{\la}_1' &=& -\widetilde{\la}_1(F'(x)-u(t) K) -
                       \eta e^{-\delta t} p u(t) K\, ,\
                       \widetilde{\la}_1(0) = \widetilde{\la}_{1,0}\, , \\
\widetilde{\la}_2' &=& \widetilde{\la}_1 x u + \gamma \widetilde{\la}_2 +
                   \eta e^{-\delta t} (c - p x) u(t)\, ,\
                   \widetilde{\la}_2(0)=\widetilde{\la}_{2,0}\, , \\
\widetilde{\la}_2(t) - \eta e^{-\delta t} r & \le & 0 \text{ for all }
                   t \in [0,\infty)\, , \\
\widetilde{\la}_2(t) - \eta e^{-\delta t} r & = & 0 \ \mu \text{--a.e. in }
                   [0,\infty)\, , \\
\widetilde{H}(t,x(t),K(t),u(t),\widetilde{\la}_1(t),\widetilde{\la}_2(t),\eta)
     & = & \max_{w\in [0,1]} \widetilde{H}(t,x(t),K(t),w,\widetilde{\la}_1(t),
           \widetilde{\la}_2(t),\eta) \text{ a.e. in } [0,\infty) \, .
\end{eqnarray*}
Now we set
$$  \la_1(t) := \widetilde{\la}_1(t) e^{\delta t}\, ,\
    \la_2(t) := \widetilde{\la}_2(t) e^{\delta t}\, ,\
    \la_{1,0}:= \la_1(0) \, ,\ \la_{2,0}:= \la_2(0) \, , $$
$$  H(t,\check{x},\check{K},w,\la,\eta) :=
    (-\la \check{x} + \eta (p \check{x}- c))\check{K} w \, . $$
Next we define the auxiliary functions $z$, $\la: [0,\infty) \to \R$
$$  z := -\la_1 x + \eta (p x - c) \, ,\ \la := \la_2 \, ,\
    z_0 := \la_{1,0} \, ,\ \la_0 := \la_{2,0} \, . $$
that are used to reinterpret the necessary conditions. We call $z$ and
$\la$ switch variables. Note that $z$ defines along the maximum condition
the value of $u(t)$, namely $u(t) = 0$ if $z(t) < 0$ and $u(t) = 1$ if
$z(t) > 0$. If $z(t)$ vanishes, then the value of $u(t)$ has to be
determined by other means. The function $\la$ defines a jump switch,
since $\mu(\{t\}) > 0$ for some $t \in [0,\infty)$ implies $\la(t) = \eta r$.

Finally we are able to rewrite the necessary optimality conditions in
the form
\begin{eqnarray*}
z_0^2 + \la_0^2 + \eta^2  &\neq&  0 \, ,\ \eta \ge 0 \, , \\[1ex]
x'  &=&  F(x)-u(t) K x \, ,\ x(0) = x_0 \, , \\[1ex]
dK  &=&  -\gamma K dt + \mu(dt) \, ,\ K(0) = K_0 \, , \\[1ex]
z'  &=&  z g(x) - \eta \psi(x) \, ,\ z(0) = z_0 \, , \\[1ex]
\la'  &=&  \kappa \la - z u(t) \, ,\ \la(0) = \la_0 \, , \\[1ex]
\la(t) - \eta r  &\le&  0 \text{ for all } t \in [0,\infty) \, , \\
\la(t) - \eta r  &=&  0 \ \mu \text{--a.e. in } [0,\infty) \, , \\[1ex]
z(t)K(t)u(t) &=& \max_{w\in [0,1]} z(t)K(t)w\ \text{ a.e. in }[0,\infty) \, .
\end{eqnarray*}
We want to exclude the irregular case $\eta = 0$. This is prepared by

\begin{lemma} \label{th:nonopt}                                      
A policy $(u,\mu)$ such that there exists $\tau > 0$ with $\mu(A) = 0$
for each measurable subset $A$ of $(\tau,\infty)$ is not optimal for each
$(x_0,K_0)\in (0,\bar{x}) \times (0,\infty)$. \\[1ex]
\bew \rm \
Suppose $(u,\mu)$ is a policy with the property that for each $\tau > 0$
we have $\mu(A)=0$ for all $A\subset [\tau,\infty)$. Then we have for the
resulting trajectory $(x,K)$: $\lim_{t \to \infty} x(t) = \bar{x}$,
$\lim_{t \to \infty} K(t) = 0$. Therefore it is enough, due to Belman's
principle of optimality, to show that such a trajectory for initial data
$(x_0, K_0)$ in a neighborhood of $(\bar{x},0)$ cannot be optimal. \\
We choose $\alpha > 0$ and $x_1 \in (x^*, \bar{x})$ with
$c_* - p x \le - \alpha$ for $x \in [x_1,\bar{x}]$; this is possible due
to the assumption (V3). Let $n \in \N$ be with $- n\alpha + p\bar{x} < 0$
and set $K_1:= F(x_1)/x_1$. \\
Let $(x,K)$ be the trajectory associated with $(u,\mu),$ i.e.
$$ x'(t) = F(x(t)) - u(t) K_0 e^{-\gamma t}x(t)\, ,\
   t > 0\, ,\ x(0) = x_0\, , $$
where $x_0 \ge x_1$, $K_0 \in (0,K_1/(n+1))$. We are able to describe a
better policy $(\breve{x}, \breve{K}, \breve{u}, \breve{\mu})$, namely:
$$ \breve{x}'(t) = F(\breve{x}(t)) - \breve{u}(t) \breve{K} \breve{x}(t)
   \, ,\ \breve{x}(0) = x_0\, ,\ \ \
   d \breve{K}(t) = -\gamma \breve{K} + \breve{\mu}(dt)\, ,\ 
   \breve{K}(0) = K_0\, , $$
where $\breve{u} \equiv 1$ and $\breve{\mu}$ describes a jump at time
$t=0$ of height $h$. Notice that we have $\breve{x}(t) \ge x_1$,
$(K_0+h) e^{-\gamma t} \le K_1$ for all $t \ge 0$. We compare the values
of the objective function: \\[2ex]
$\D rh + \int_0^\infty e^{-\delta t} (c-p\breve{x}(t)) \breve{u}(t)
 \breve{K}(t)
 \, dt - \int_0^\infty e^{-\delta t} (c-px(t)) u(t) K_0 e^{-\gamma t}\, dt$
\begin{eqnarray*}
&\le& rh + \int_0^\infty e^{-\delta t}(c-p\breve{x}(t)) (K_0+h)
      e^{-\gamma t}\, dt
      - \int_0^\infty e^{-\delta t} (c-px(t)) K_0 e^{-\gamma t}\, dt \\
&=& rh + \int_0^\infty e^{-(\delta + \gamma)t} (c-p\breve{x}(t)) (K_0+h)\, dt
    - \int_0^\infty e^{-(\delta + \gamma)t} (c-p\breve{x}(t)) K_0\, dt \\
& & + \int_0^\infty e^{-(\delta + \gamma)t} (c-p\breve{x}(t))K_0\, dt
    - \int_0^\infty e^{-(\delta + \gamma)t} (c-px(t)) K_0\, dt \\
&=& h \int_0^\infty e^{-(\delta + \gamma)t} r'\, dt 
    + h \int_0^\infty e^{-(\delta + \gamma)t} (c-p\breve{x}(t))\, dt 
    + K_0 p \int_0^\infty e^{-(\delta + \gamma)t} (x(t)-\breve{x}(t))\, dt \\
&=& h \int_0^\infty e^{-(\delta + \gamma)t} (c_* - p\breve{x}(t))\, dt + K_0
    p \int_0^\infty e^{-(\delta + \gamma)t} (x(t)-\breve{x}(t))\, dt \\[1ex]
&\le& -h\alpha (\delta+\gamma)^{-1} + K_0 p \bar{x} (\delta+\gamma)^{-1}
\ = \ K_0 (\delta+\gamma)^{-1} (-n \alpha + p \bar{x})
\ < \ 0\, .
\end{eqnarray*}
\eop
\end{lemma}

\begin{folge} \label{kor:nonopt0}                                    
Let $(x,K,u,\mu)$ be an optimal process with adjoint variables $\eta$, $z$,
$\la$. If $\eta = 0$, then there exists for any $\tau > 0$ some $t > \tau$
with $\la(t) = 0$. \\[1ex]
\bew \rm \
If there exists a $\tau >0$ such that $\la(t) < 0$, for all $t > \tau$,
then $\mu(A) = 0$ for each measurable subset $A$ of $(\tau,\infty)$. This
contradicts the optimality of the policy $(u,\mu)$ (see Lemma~%
\ref{th:nonopt}), proving the corollary.
\eop
\end{folge}

\begin{satz} \label{th:eta0}                                       
Let $(x,K,u,\mu)$ be an optimal process with adjoint variables $\eta$, $z$,
$\la$. Then $\eta \neq 0$. \\[1ex]
\bew \rm \
If $\eta = 0$, we have the necessary conditions $z_0^2 + \la_0^2 \neq 0$,
\begin{eqnarray*}
z'   &=& z g(x) \, ,\ z(0) = z_0 \ , \\
\la' &=& \kappa \la - z u(t) \, ,\ \la(0) = \la_0 \, , \\
\la(t) &\le& 0 \text{ for all } t \in [0,\infty) \, , \\
\la(t) & = & 0 \ \mu \text{--a.e. in } [0,\infty) \, .
\end{eqnarray*}
We distinguish six cases according to the initial conditions $(\la_0,z_0)$: \\
{\it i)} \ $\la_0 = 0$, $z_0 = 0$: \ Here we have $z(t) = 0$, $\la(t) = 0$,
for all $t > 0$. This contradicts the necessary condition $z_0^2 + \la_0^2
\neq 0$. \\
{\it ii)} \ $\la_0 = 0$, $z_0 < 0$: \ Clearly, $\la'(0) = 0$ and $z'(0) < 0$.
This implies $z(t)<0$, $u(t)=0$, $\la(t)=0$ for all $t>0$. Therefore, the
policy $(u,\mu)$ is not better than $(u,\hat{\mu})$ with $\hat{\mu} \equiv 0$.
 From Corollary~\ref{kor:nonopt0} we know that already this policy is
not optimal. \\
{\it iii)} \ $\la_0 = 0$, $z_0 > 0$: \ Note that $z'(0) > 0$, $\la'(0) < 0$
and we have $z(t) > 0$, $u(t) = 1$, $\la(t) < 0$ for all $t > 0$. From
Corollary~\ref{kor:nonopt0} we obtain that this policy is not optimal. \\
{\it iv)} \ $\la_0 < 0$, $z_0 = 0$: \ We have $z'(0) = 0$, $\la'(0) < 0$.
This implies $z(t) = 0$, $\la(t) < 0$ for all $t > 0$. From Corollary~%
\ref{kor:nonopt0} follows that this policy is not optimal. \\
{\it v)} \ $\la_0 < 0$, $z_0 < 0$: \ We have $z'(0) < 0$, $\la'(0) < 0$.
This implies $z(t) < 0$, $\la(t) < 0$ for all $t > 0$ and again due to
Corollary~\ref{kor:nonopt0} this policy is not optimal. \\
{\it vi)} \ $\la_0 < 0$, $z_0 > 0$: \ In this case $z'(0) > 0$, $\la'(0) < 0$
and we have $z(t) > 0$, $\la(t) < 0$ for all $t > 0$. Applying Corollary~%
\ref{kor:nonopt0}, we conclude that this policy is not optimal. \\
Therefore, in all cases we have a contradiction and the theorem is proved.
\eop
\end{satz}
%
%
%
\section{Exploitation of the necessary conditions} \label{sec:exploit}

An immediate consequence of Theorem~\ref{th:eta0} is the fact that the
Lagrange multiplier $\eta$ can be chosen equal to one. In this case we
have to analyze the following necessary conditions:
\begin{eqnarray*}
x' &=& F(x) - u(t) K x \, ,\ x(0) = x_0\, , \\
dK &=& -\gamma K dt + \mu(dt) \, ,\ K(0) = K_0 \, , \\
z' &=& z g(x) - \psi(x) \, ,\ z(0) = z_0 \, , \\
\la' &=& (\la -r) \kappa - z  u(t) + r' \, ,\ \la(0) = \la_0 \, , \\
\la(t) &\le& r \, , \mbox{ for all } t \in [0,\infty) \, , \\
\la(t) & = & r \, , \ \mu \mbox{--a.e. in } [0,\infty) \, , \\
z(t)K(t)u(t) &=& \max_{w\in [0,1]} z(t)K(t)w \ \mbox{ a.e. in } [0,\infty)\, .
\end{eqnarray*}
The central problem in the analysis of the necessary conditions consists in
finding out the initial conditions $\la_0 = \la(0)$, $z_0 = z(0)$ which are
in agreement with the condition \\[2ex]
$(R)$ \hfil $\la(t) \le r$, for all $t \in [0,\infty)$. \\[2ex]
One can easily check that if $\la(\tau) = r$ for some $\tau > 0$, then
$\la'(\tau) = 0$ and $\la''(\tau) \le 0$. It is also clear that $K_{0,+}
\neq K_0$ at most if we have $\la(0) = r$. Knowing that it is possible to
chose $\eta = 1$, we formulate Corollary~\ref{kor:nonopt0} again:

\begin{folge} \label{kor:nonopt1}                                    
Let $(x,K,u,\mu)$ be an optimal process with adjoint variables $z$, $\la$.
Then for each $\tau > 0$ there exists $t > \tau$ with $\la(t) = r$. \\[1ex]
\bew \rm \
See the proof for Corollary~\ref{kor:nonopt0}.
\eop
\end{folge}

The rest of this section is devoted to the analysis of the relationship
between the initial states $(x_0,K_0)$ and the initial conditions $(z_0,
\la_0)$ of the adjoint variables. We prove a series of auxiliary lemmas,
that will allow us in the next section to detect the optimal trajectories
and correspondent policies.

\begin{lemma} \label{lem:sprung1}                                    
Let $(x,K,u,\mu)$ be optimal with adjoint variables $z$, $\la$. Then there
exists no $\tau > 0$ with $x(\tau) > x^*$ and $F(x(\tau)) \ge K(\tau)
x(\tau)$. \\[1ex]
\bew \rm \
Assume the contrary. We consider two cases separately: \\
{\it i)} \ If $\la(\tau) = r$, then it follows from $(R)$ that
$\la'(\tau) = 0$. Consequently, $z(\tau) u(\tau) = r'$ and $z(\tau) = r'$.
Now, it follows from $(V4)$ and the assumption $x(\tau) > x^*$ that
$\psi_*(x(\tau)) > 0$ and substituting in the dynamic of $z$ we have
$$ z'(\tau) = (z(\tau)-r') g(x(\tau)) - \psi_*(x(\tau)) < 0 \, . $$
Because of $\la''(\tau) = -z'(\tau) > 0$ we have a contradiction to $(R)$. \\
{\it ii)} \ If $\la(\tau) < r$, then we obtain from Corollary~%
\ref{kor:nonopt1} a $\tau_1 > \tau$ with $\la(\tau_1) = r$, $\la(t) < r$,
for $t \in [\tau, \tau_1)$. Observing the dynamic of the pair $(x(t),K(t))$,
we conclude that $0 \le F(x(t)) - K(t)x(t) \le x'(t)$, $t \ge \tau$. Then
for $t=\tau_1$ we have $F(x(\tau_1)) \ge K(\tau_1) x(\tau_1)$, $x(\tau_1)
\ge x(\tau) > x^*$. To obtain a contradiction, we argument as in {\it i)},
setting $\tau := \tau_1$.
\eop
\end{lemma}

\begin{lemma} \label{lem:lrplus}                                     
Let $(x,K,u,\mu)$ be optimal with adjoint variables $z$, $\la$. Then for
$x_0 = x^*$, $K_0 > K^*$ the initial conditions $\la_0 = r$, $z_0 \ge r'$
are not possible. \\[1ex]
\bew \rm \
Assume that $\la_0 = r$, $z_0 = r'$. From the differential equations for
$x$, $\la$ and $z$ follows $x'(0) < 0$, $\la'(0) = 0$, $z'(0) = 0$. Then
we have
$$ \la''(0)   = \la'(0) \kappa - z'(0) = 0 \, ,\ \
   \la'''(0+) = \psi_*'(x^*) [F(x^*) - K_{0,+} x^*] < 0 \, . $$
Therefore, $\la(t) \not\equiv r$ does not occur. This holds also if $z_0
> r'$ since in this case we have $\la'(0) < 0$. \\
Due to Corollary~\ref{kor:nonopt1}, $\la(t) < r$ for all $t > 0$ is not
possible. Let
$$ \tau := \sup\{ \sigma \ge 0\, |\, \la(t) < r,\ 0<t<\sigma \} \, . $$
We know that $\tau < \infty$, $\la(\tau) = r$, $\la'(\tau) = 0$,
$z(\tau) = r'$. Since $x_0 = x^*$ and $x'(0) < 0$, we have three
possible cases: \\
{\it i)} \ If $x(t) < x^*$ for all $t \in (0,\tau]$, then $z'(t) > 0$, for
$t \in [0, \tau)$ and $z(0) = z(\tau) = r'$ which is not possible. \\
{\it ii)} \ If $x(t) < x^*$, for all $t \in (0,\tau)$ and $x(\tau) = x^*$,
then we may argument as in {\it i)}. \\
{\it iii)} \ If $x(\tau) > x^*$, then we can choose $\sigma \in (0,\tau)$
with $x(t) < x^*$, for $t \in [0, \sigma)$, and $x(\sigma) = x^*$. Note that
$z(t) > r'$ (and consequently $u(t) = 1$) for $t \in (0,\sigma]$. From the
definition of $\sigma$ we have $x'(\sigma) \ge 0$, i.e. $K(\sigma) \le
F(x^*)/x^*$. Now, arguing as in case {\it ii)} in the proof of Lemma~%
\ref{lem:sprung1}, we obtain the inequality $K(\sigma) > F(x^*)/x^*$,
again a contradiction.%
\footnote{One should note that the development in case {\it ii)} in the 
proof of Lemma~\ref{lem:sprung1} still holds if $x(\tau) = x^*$.}
\eop
\end{lemma}

\begin{lemma}\label{lem:lr}                                         
Let $(x,K,u,\mu)$ be optimal with adjoint variables $z$, $\la$. If
$x_0 = x^*$, $K_0 = K^*$, then $\mu = \gamma K^* dt$ and $x(t) = x^*$,
$K(t) = K^*$, $\la(t) = r$, $z(t) = r'$, $u(t) = 1$, for $t \ge 0$. \\[1ex]
\bew \rm \
Assume $\la_0 < r$. From Corollary~\ref{kor:nonopt1} we obtain $\tau > 0$
with $\la(\tau) = r$, $\la(t) < r$, for $t \in (0, \tau)$. Condition (R)
implies $\la'(\tau) = 0$ and therefore $z(\tau) = r'$; especially we have
$u(t) \equiv 1$ in a neighborhood of $t = \tau$. From $x' = F(x) - K x$ we
obtain $x(\tau) > x^*$ and $\la''(\tau) = -z'(\tau) = \psi_*(x(\tau)) > 0$,
contradicting (R). Therefore, we must have $\la_0 = r$. Next verify that
$z_0 = r'$: \\
{\it i)} \ If $z_0 < r'$, then $\la'(0) = r' - z_0 > 0$, contradicting (R). \\
{\it ii)} \ If $z_0 > r'$, then $\la'(0) < 0$. Arguing with Corollary~%
\ref{kor:nonopt1} (see the beginning of this proof), we obtain a $\tau > 0$
with $\la(\tau) = r$, $\la'(\tau) = 0$ but $\la''(\tau) > 0$, again
contradicting (R). \\
Therefore, we must have $z_0 = r'$. Now define
$$ \tau := \sup\{ \sigma \ge 0\, |\, \la(t)=r,\ 0 \le t \le \sigma \} \, . $$
If $\tau = 0$, then $\la'(0) < 0$ and arguing as in {\it ii)} above we
obtain a contradiction. Therefore, we must have $\tau > 0$. Note that
if $\sigma \in (0,\tau)$, then $\la(t) = r$, $\la'(t) = 0$, $z(t) = r'$,
$z'(t) = 0$, for $t \in [0,\sigma]$. It follows $\psi_*(x(t)) = 0$, for
$t \in [0,\sigma]$, and with condition (V4) we have $x(t) = x^*$, for
$t \in [0,\sigma]$. From the differential equation $x' = F(x)- K x$
follows $K(t) x^* = F(x^*)$, for $t \in [0,\sigma]$, and with $dK =
-\gamma K dt + \mu(dt)$ we finally obtain $\mu_{|_{[0,\sigma]}} =
-\gamma K^* dt$. \\
Clearly, it is enough to prove $\tau = \infty$. If this were not the case,
we would have either $K(\tau+) > K^*$ or $K(\tau+) = K^*$. Due to Lemma~%
\ref{lem:lrplus}, $K(\tau+) > K^*$ is not possible. The other case,
$K(\tau+) = K^*$, cannot be true, since otherwise we could repeat the
complete argumentation for the initial time $t=\tau$, contradicting the
maximality of $\tau$.
\eop
\end{lemma}

\begin{lemma} \label{lem:lrminus}                                   
Let $(x,K,u,\mu)$ be optimal with adjoint variables $z$, $\la$. Then for
$x_0 = x^*$, $K_0 \le K^*$ we must have the initial values $\la_0 = r$,
$z_0 = r'$, $\la'(0) = 0$, $z'(0) = 0$. Further we have $K_{0,+} = K^*$.
\\[1ex]
\bew \rm \
The equality $\la'(0) = 0$, $z'(0) = 0$ follow from $x_0 = x^*$, $\la_0 = r$,
$z_0 = r'$ in an obvious way. Consequently, we only have to prove $\la_0 = r$,
$z_0 = r'$. \\
If $\la_0 < r$, then there exists $\tau > 0$ such that $\la(t) < r$, for
$t \in [0,\tau)$. Then $K(t) < K_0 \le K^*$ for $t \in (0,\tau]$, and
consequently $x(\tau) > x^*$, $x'(\tau) = F(x(\tau)) - K(\tau) x(\tau)
> 0$. But this cannot occur due to Lemma~\ref{lem:sprung1}. Therefore,
$\la_0 = r$. \\
If $z_0 < r'$, then $\la'(0) < 0$. Arguing as before (assumption $\la_0
< r$) we obtain a contradiction. Therefore, $z_0 \ge r'$. Next we exclude
the case $z_0 > r'$. \\
If $K_0 = K^*$, then Lemma~\ref{lem:lr} implies $z_0 = r'$ proving the
theorem. If $K_0 < K^*$, we consider the following cases: \\
{\it i)} \ $K_{0,+} > K^*$ is not possible due to Lemma~\ref{lem:lrplus}. \\
{\it ii)} \ If $K_{0,+} = K^*$, then $z_0 = r'$ follows from Lemma~%
\ref{lem:lr}. \\
{\it iii)} \ $K_{0,+} \in [K_0,K^*)$. We obtain from the dynamic of the pair
$(x,K)$ some $\tau > 0$ with $x(\tau) > x^*$, $F(x(\tau)) - K(\tau) x(\tau)
> 0$. But this is not possible due to Lemma~\ref{lem:sprung1}. \\
Therefore, only $K_{0,+} = K^*$ is possible and the theorem is proved.
\eop
\end{lemma}

\begin{lemma} \label{lem:sprung2}                                   
Let $(x,K,u,\mu)$ be optimal with adjoint variables $z$, $\la$. If $x_0 >
x^*$ and $K_0 < F(x_0) / x_0$, then we must have: $\la_0 = r$, $z_0 > r'$,
$K_{0,+} > K_0$. \\[1ex]
\bew \rm \
Since $K_{0,+} = K_0$ is not allowed (see Lemma~\ref{lem:sprung1}), we
have $K_{0,+} > K_0$ and $\la_0 = r$. Then $\la'(0) \le 0$, and we obtain
$z_0 \ge r'$. If $z_0 = r'$, we would have $z'(0) = - \psi_*(x_0) < 0$,
$\la'(0) = 0$ and $\la''(0) = -z'(0) > 0$. This however contradicts the
condition (R).
\eop
\end{lemma}

\begin{figure}[!h]
\begin{center}
\epsfig{figure=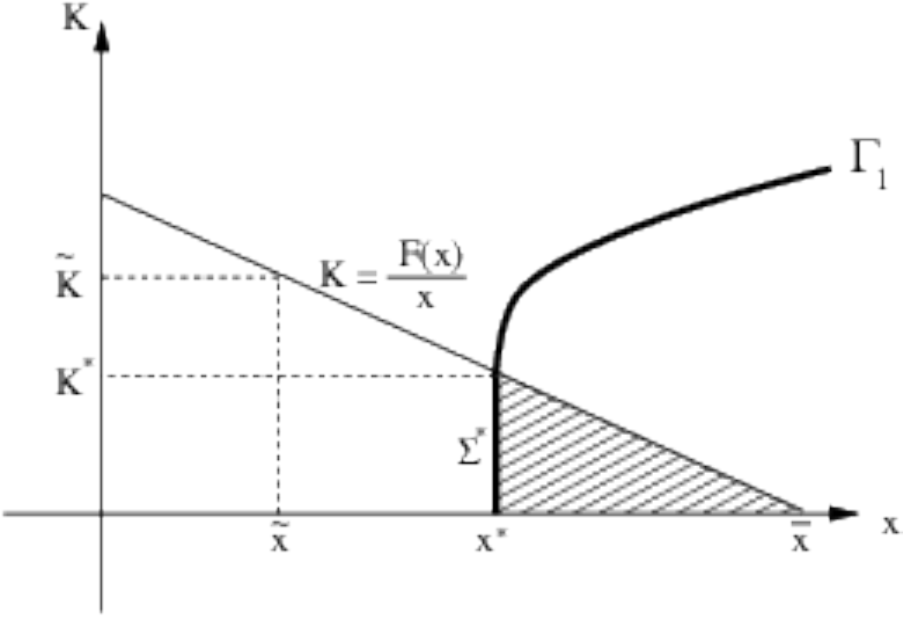,width=.75\textwidth}
\caption{Curves $\Gamma_1$ and $\Sigma^*$.} \label{fig:gamma1}
\end{center}
\end{figure}

Let $(x,K)$ be a solution of the system
\begin{equation} \label{eq:gam1}
\left\{ \begin{array}{l}
  x' = -F(x) + Kx, \ x(0) = x^* \\
  K' = \gamma K, \ K(0) = K^* \, ,
\end{array} \right.
\end{equation}
with interval of existence $[0,\tau)$. Since $K(t) = K^* e^{\gamma t}$,
$x'(0) = 0$ and $x''(0) = K^* e^{\gamma t} \gamma x^* > 0$, we have
$(x(t),K(t)) \in (x^*,\bar{x}) \times (K^*,\infty)$ for $t>0$ small.
Therefore, it is easy to see that there exists some $\bar{t} \in (0,\tau)$
with $x(\bar{t}) = \bar{x}$ and $x^* \le x(t) \le \bar{x}$, $t \in
[0,\bar{t}]$. The curve defined by
$$ [0,\bar{t}] \ni t \longmapsto
   (x(t),K(t)) \in [x^*,\bar{x}] \times [K^*,\infty) $$
is denoted by $\Gamma_1$. Note that at $t = 0$ we have $(x'(0),K'(0)) =
(0,K^* \gamma)$, where $K^* \gamma > 0$. In Figure~\ref{fig:gamma1} we
illustrate the construction of the curve $\Gamma_1$ for the case of the
logistic function.

\begin{lemma} \label{lem:skal}                                      
There exists a function $h_1: [x^*,\bar{x}] \to [K^*,\infty)$ such that
\begin{itemize}
\item [(a)] $\Gamma_1 = \{ (x,h_1(x))\ |\ x \in [x^*,\bar{x}] \}$;
\item [(b)] $h_1$ is twice continuous differentiable in $(x^*, \bar{x})$ 
            and monotone increasing.
\end{itemize}
\bew \rm \
Note that the system \eqref{eq:gam1} can be transformed into the scalar
equation:
\begin{equation} \label{eq:gam2}
\frac{dK}{dx} = \frac{\gamma K}{-F(x) + Kx} \, ,\ K(x^*) = K^* \, .
\end{equation}
It becomes clear that $\Gamma_1$ has a parameterization $[x^*,\bar{x}]
\ni x \mapsto (x, h_1(x)) \in [x^*,\bar{x}] \times [K^*,\infty)$ and the
assertions follow.
\eop
\end{lemma}

\begin{lemma} \label{lem:sprung3}                                   
Let $(x,K,u,\mu)$ be optimal with adjoint variables $z$, $\la$. For
$x_0 > x^*$ and $K_0 < h_1(x_0)$, the initial conditions have to satisfy
$\la_0 = r$, $z_0 > r'$, $K_{0,+} \ge h_1(x_0)$. \\[1ex]
\bew \rm \
Assume $\la_0 < r$.
We define $\tau := \sup\{ \sigma \ge 0\ |\ \la(t) < r,\ t \in [0,\sigma) \}$.
Due to Corollary~\ref{kor:nonopt1}, $\tau$ is positive and finite.
Then we have $\la(\tau) = r$, $\la'(\tau) = 0$, $z(\tau) = r'$.
Considering the definition of $\Gamma_1$, we obtain $x(t) > x^*$,
$t \in (0,\tau]$ due to the values of $\la$.
Now, from assumption (V4), the differential equation for $z$ and
$x(\tau) > x^*$, follow $z'(\tau) < 0$.
Again from the differential equation for $z$ we obtain $z(t) > r'$,
$t \in (0,\tau)$.
This implies $\la'(t) < 0$, $t \in (0,\tau)$, which is a contradiction to
the definition of $\tau$. Therefore, $\la_0 = r$ must occur.
This implies $0 \ge \la'(0) = r' - z_0$ and $z_0 \ge r'$ follows. \\
If $z_0 = r'$, we would have $z'(0) = -\psi_*(x_0) < 0$, $\la'(0) = 0$ and
$\la''(0) = -z'(0) > 0$. This however contradicts the condition (R).
Therefore, we must have $z_0 > r'$. \\
Assume $K_{0,+} = K_0$. We know already that $z_0 > r'$. Then $\la'(0) < 0$
holds and by the same arguments as above (see $\la_0 < r$) we obtain a
contradiction. Therefore, $K_{0,+} > K_0$ must hold. If $K_{0,+} \in
(x_0, h_1(x_0))$ we repeat the argumentation above with $K_0 :=
K_{0,+} < h_1(x_0)$, obtaining again a contradiction. Thus, we must
have $K_{0,+} \ge h_1(x_0)$.
\eop
\end{lemma}

We denote the curve
$$ [0,K^*] \ni K \longmapsto (x^*,K) \in [0,\infty) \times [0,K^*] $$
by $\Sigma^*$. Next we verify that when an optimal trajectory meets the
curve $\Sigma^*$, some properties have to be satisfied.

\begin{lemma} \label{lem:erreich}                                   
Let $(x,K,u,\mu)$ be optimal with adjoint variables $z$, $\la$. Let
$x(\sigma) = x^*$, $K(\sigma) \in (0,K^*)$ for some $\sigma > 0$. Then
$$ \la(\sigma) = r \, ,\ \la'(\sigma)= 0 \, ,\
   z(\sigma) = r' \, ,\ z'(\sigma)=0 \, . $$
\bew \rm \
Assume $\la(\sigma)< r$. From the differential equation for $x$ and $K$
we obtain $\tau > \sigma$ with
$$ x(\tau) > x^* \, ,\ F(x(\tau)) - K(\tau)x(\tau) > 0 \, , $$
which is in contradiction to Lemma~\ref{lem:sprung1}. Thus we must have
$\la(\sigma) = r$. Therefore, $\la'(\sigma) = 0$, $z(\sigma) = r'$. Finally,
$z'(\sigma) = 0$ follows from $\psi_*(x^*) = 0$.
\eop
\end{lemma}

In the neighborhood of $\Gamma_1$ and $\Sigma^*$ we have obtained a lot
of information concerning the behavior of an extremal trajectory. Lemma~%
\ref{lem:sprung3}, ensures that if $x_0 > x^*$ and $K_0 < h_1(x_0)$, there
must be a jump at $t = 0$. Furthermore Lemma~\ref{lem:sprung1} says that an
optimal trajectory $(x(t), K(t))$ does not enter the dashed region in
Figure~\ref{fig:gamma1}.

Our next step is to analyze the behavior of the optimal trajectories that
meet the curve $\Sigma^*$, i.e. $(x(\tau),K(\tau)) \in \Sigma^*$ for some
$\tau > 0$. Since the curve $\Sigma^*$ is reached by an optimal trajectory
with $z = r'$, $\la = r$, $x = x^*$ and $K \in [0,K^*]$, we use this
information to solve the differential equations for $x$, $K$, $z$ and
$\la$ backwards in time:
$$ \left\{ \begin{array}{rl@{}l}
     x'  & = -F(x) + K x                & ,\  x(0) = x^* \, ,\\
     K'  & = \gamma K                   & ,\  K(0) = K_1 \, , \\
     z'  & = -(z-r') g(x) + \psi_*(x)   & ,\  z(0) = r' \, ,\\
     \la' & = -(\la -r) \kappa + z - r' & ,\ \la(0) = r \, ,
   \end{array} \right. $$
where $K_1 \in [0,K^*]$. Such a trajectory eventually comes close to
$(x,K) = (0,0)$ where we expect an optimal control $u \equiv 0$.
Therefore, the zeros of the switching variable $z$ are of interest in
this region. Since for $z(t) < r'$ a value $\la(t) = r$ is not
allowed, the behavior of the adjoint variable $\la$ is therefore not
so important in this region.

\begin{lemma} \label{lem:nullst}                                    
For each $K_1 \in [0,K^*]$, let $(x,z)$ be the solution of
\begin{equation}\label{eq:nullst}
\left\{ \begin{array}{l@{}l}
  x' = -F(x)+ K_1e^{\gamma t} x & ,\ x(0) = x^* \, , \\
  z' = -(z-r') g(x) + \psi_*(x) & ,\ z(0) = r' \, .
\end{array} \right.
\end{equation}
Then there exists for each $K_1 \in (0,K^*)$  $\tau := \tau_{K_1} > 0$
with $x(\tau) = x^*$. Moreover there exists $\widetilde{K}_1 \in (0,K^*)$
such that the following assertions hold:
\begin{itemize}
\item [(a)] If $K_1 \in (0,\widetilde{K}_1)$, then $z$ has a uniquely
            determined zero $\sigma \in (0,\tau)$, where $z'(\sigma) < 0$
            and $x(\sigma) \in (0,\tilde{x})$;
\item [(b)] If $K_1 \in (\widetilde{K}_1,K^*)$, then $z(t) > 0$ for all
            $t \in [0,\tau]$;
\item [(c)] If $K_1 = \widetilde{K}_1$, then there is a uniquely determined
            $\sigma$ in $(0,\tau)$ with $z(\sigma) = z' (\sigma) = 0$;
            moreover $x(\sigma) = \tilde{x}$ and $\ttilde{K} :=
            \widetilde{K}_1 e^{\gamma \sigma} > F(\tilde{x}) / \tilde{x}
            = \widetilde{K}$.
\end{itemize}
\bew \rm \
Consider the solution $(x,z)$ of \eqref{eq:nullst} with $K_1 = 0$. Since
$x_0 = 0$ is an attracting equilibrium point of $x' = -F(x)$ (see Remark~%
\ref{rem:red1}), the solution $x$ exists for all times $t \ge 0$ and
$\lim_{t\to\infty} x(t) = 0$. Due to this fact $z$ is also defined for
all $t \ge 0$. Assume $z(t) > 0$ for all $t >0$. As we know, the
differential equation may be formulated as $z' = -z g(x) + \psi(x)$.
Since $\psi$ is negative and continuous in $[0,\tilde{x})$ (see (V4)),
there exists some $a > 0$ such that
$$ \psi(\xi) \le -a,\ {\rm for\ } \xi \in [0,\T\frac{\tilde{x}}{2}] \, . $$
Since $\D\lim_{t\to\infty} x(t) = 0$, there is some $t_0 > 0$ with $x(t)
\in [0,\frac{\tilde{x}}{2}]$, $t \ge t_0$. This implies for $t \ge t_0$
$$ z(t) - z(t_0)  \ = \ \int_{t_0}^t \big[ -z(s)g(x(s)) + \psi(x(s)) \big] ds
                \ \le \ \int_{t_0}^t \psi(x(s)) ds \ \le \ -a(t-t_0) \, . $$
But this contradicts the hypothesis on $z$. Thus, there must be a
$\sigma > 0$ with $z(\sigma) = 0$ and $z(t) > 0$, $t \in [0,\sigma)$.
Then $z'(\sigma) \le 0$, $\psi(x(\sigma)) \le 0$, and therefore
$x(\sigma) \le \tilde{x}$. Due to
$$ 0 \le z''(\sigma) = -z'(\sigma) g(x(\sigma)) + \psi'(x(\sigma)) x'(\sigma)
      = -z'(\sigma) g(x(\sigma)) - \psi'(x(\sigma)) F(x(\sigma)) $$
and assumption (V5), $z'(\sigma) = 0$ cannot occur. Therefore, $z'(\sigma)
< 0$, $x(\sigma) < \tilde{x}$ and $z(t) < 0$ for $t > \sigma$ due to the
differential equation for $z$. By continuity arguments, we may choose
a maximal $\widetilde{K}_1 > 0$ such that for the solution $(x,z)$ of
\eqref{eq:nullst} with $K_1 \in [0,\widetilde{K}_1)$ there exists
$\tau > 0$ and $\sigma \in (0,\tau]$ with $x(\tau) = x^*$,
$$ z(t) > 0 \, ,\ t \in [0,\sigma) \, ,\ \
   z(\sigma) = 0 \, ,\ \
   z(t) < 0 \, ,\ t \in (\sigma,\tau] \, ; $$
proving {\it (a)} and {\it (b)}.\\
Now we prove  {\it (c)}. From the construction of $\widetilde{K}_1$
and due to the differential equation for $x$ we obtain $z(\sigma) = 0$ and
$x(\sigma) > 0$. Since $\widetilde{K}_1$ is maximal, we have $ z(\sigma) =
z'(\sigma) = 0$. We cannot have $x(\sigma) > \tilde{x}$, since $z'(\sigma)
= - \psi(x(\sigma)) < 0$. The case $x(\sigma) < \tilde{x}$ cannot be
occur, since $z'(\sigma) = -\psi(x(\sigma)) > 0$. Therefore, we must have
$x(\sigma) = \tilde{x}$. \\
We cannot have $K(\sigma) x(\sigma) < F(x(\sigma))$ since, due to
$z''(\sigma) = \psi'(x(\sigma)) x'(\sigma) < 0$, $z$ would be negative
in a neighborhood of $\sigma$. Otherwise, if we had $K(\sigma) x(\sigma)
= F(x(\sigma))$, i.e. $K(\sigma) \tilde{x} = F(\tilde{x})$, then $z$
would be positive in a neighborhood of $\sigma$, due to $z''(\sigma) = 0$
 and $z'''(\sigma) = \psi'(x(\sigma)) x''(\sigma) = \psi'(x(\sigma))
\gamma F(\tilde{x}) > 0$. Thus, we must have $K(\sigma) > F(\tilde{x}) /
\tilde{x}$ and {\it (c)} is proved.
\eop
\end{lemma}

\begin{folge} \label{kor:nullkurve}                                 
Let $\widetilde{K}_1 \in (0,K^*)$ be chosen as in Lemma~\ref{lem:nullst}.
Let $K_1 \in [0,\widetilde{K}_1]$ and $(x,z)$ the corresponding solution
of \eqref{eq:nullst}. Then there exists a continuous mapping $l:
[0,\widetilde{K}_1] \to (0,\infty)$ and $\breve{x} \in (0,\tilde{x})$,
such that
\begin{itemize}
\item [(a)] $x(l(0)) = \breve{x}$, $x(l(\widetilde{K}_1)) = \tilde{x}$;
\item [(b)] $z(l(K_1)) = 0$, $K_1 \in [0,\widetilde{K}_1]$, \ \
            $z'(l(K_1)) < 0$, $K_1 \in [0,\widetilde{K}_1)$, \ \
            $z'(l(\widetilde{K}_1)) = 0$;
\item [(c)] $l$ is continuous differentiable in $(0,\widetilde{K}_1)$
            and $l'(K_1) > 0$, $K_1 \in (0,\widetilde{K}_1)$.
\end{itemize}
\bew \rm \
Let $K_1 = 0$ and let $\sigma > 0$ with $z(\sigma) = 0$ (see proof of
Lemma~\ref{lem:nullst}). Now define $\breve{x} := x(\sigma)$. Given
$K_1 \in (0,\widetilde{K}_1)$, we denote by $x := x(\cdot;K_1)$,
$z := z(\cdot;K_1)$ the corresponding solution of \eqref{eq:nullst}
and define $l(K_1) := \sigma$, where $\sigma > 0$ is the uniquely
determined zero of $z(\cdot;K_1)$ (see Lemma~\ref{lem:nullst}); then
{\it (b)} holds. The mapping $K_1 \mapsto l(K_1)$ is continuous since
$z$ depends continuously on $K_1$ (notice that $z'(\sigma) < 0$). Now we
may extend $l: (0,\widetilde{K}_1) \to (0,\infty)$ to a continuous map on
$[0,\widetilde{K}_1]$, proving {\it (a)}. Notice that the mapping $\Psi:
[0,\tau] \times [0,\widetilde{K}_1] \to \R$ defined by $\Psi(s,K_1) :=
z(s;K_1)$ is differentiable and
$$ \frac{\partial \Psi}{\partial s}(s,K_1) = z'(s;K_1) \, ,\ \
   \frac{\partial \Psi}{\partial K_1}(s,K_1) = v(s;K_1) \, , $$
where $v:= v(\cdot;K_1):=\frac{\D\partial z}{\D\partial K_1}(\cdot;K_1)$
is the solution of the initial value problem
$$  v' = - v g(x(t;K_1)) + \underset{(*)} {\underbrace{ \big[
         -(z(t;K_1)-r') g'(x(t;K_1)) + \psi_*'(x(t;K_1)) \big] } } w \, ,
\,v(0)=0 \, , $$
and $w := \frac{\D\partial x}{\D\partial K_1}$ solves
$$ w' = \big[ -F'(x(t;K_1)) + K_1 e^{\gamma t} x(t;K_1) \big] w +
          e^{\gamma t} x(t;K_1) \, , \,w(0)=0 \, . $$
By the implicit function theorem applied on $\Psi(l(K_1),K_1) = 0$
for $K_1 \in (0,\widetilde{K}_1)$ we have
\begin{equation} \label{eq:zeq}
z'(l(K_1);K_1) \, l'(K_1) = -v(l(K_1);K_1) \, ,\ \
z'(l(K_1);K_1) < 0 \, ,\ K_1 \in (0,\widetilde{K}_1) \, ,
\end{equation}
It is obvious that $v(t) > 0$ for all $t > 0$, since $(*)$ as well as
$w(t)$ are positive for $t \in (0,l(K_1))$. From \eqref{eq:zeq} we obtain
$l'(K_1) > 0$, $K_1 \in (0,\widetilde{K}_1)$, and {\it (c)} is proved.
\eop
\end{folge}

\begin{figure}[!h]
\begin{center}
\epsfig{figure=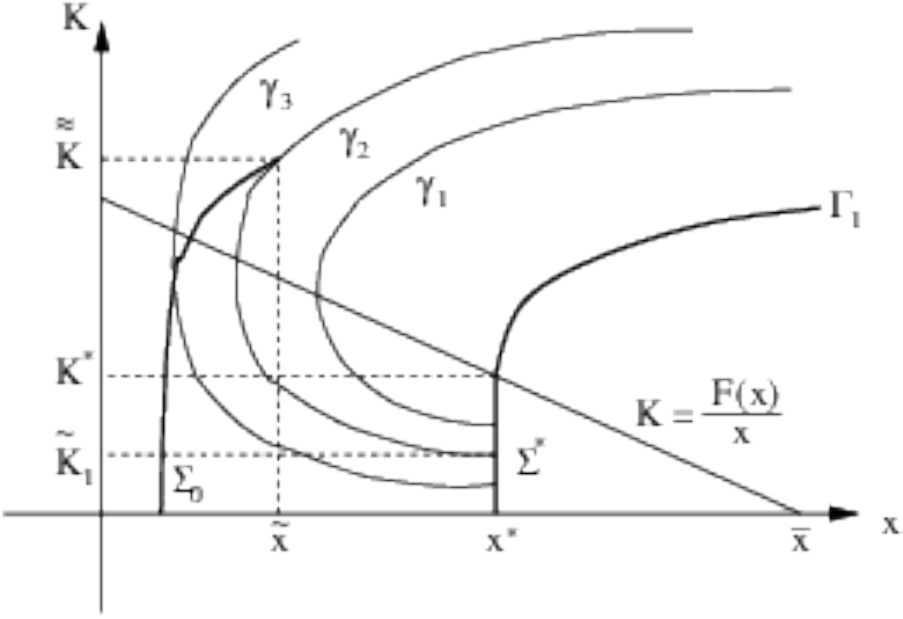,width=.75\textwidth}
\caption{Curve $\Sigma_0$.} \label{fig:sigma0}
\end{center}
\end{figure}

Corollary~\ref{kor:nullkurve} allow us to define a curve $\Sigma_0$,
parameterized by
\begin{eqnarray*}
\xi: [0,\widetilde{K}_1] & \longrightarrow &
     [0,\tilde{x}] \times [0,\ttilde{K}] \\
K_1 & \longmapsto & (x(l(K_1);K_1),K_1 e^{\gamma l(K_1)})
\end{eqnarray*}
(see the proof of Corollary~\ref{kor:nullkurve} for the notation). From
Lemma~\ref{lem:nullst} and Corollary~\ref{kor:nullkurve} we conclude that
for the solutions $(x,K,z)$ of
$$ \left\{ \begin{array}{l@{}l}
     x' = -F(x) + K x              & ,\ x(0) = x^* \, , \\
     z' = -(z-r') g(x) + \psi_*(x) & ,\ z(0) = r' \, , \\
     K'  = \gamma K              & ,\ K(0) = K_1 \in [0,K^*]
   \end{array} \right. $$
one of the following alternatives holds: \\
{\it (i)} if $K_1 \in (\widetilde{K}_1, K^*]$, then $z(t) > 0$, for all
$t \ge 0$ (see curve $\gamma_1$ in Figure~\ref{fig:sigma0}); \\
{\it (ii)} if $K_1 = \widetilde{K}_1$, then $z(t) > 0$ except at a single
time point, where $(x,K) = (\tilde{x}, \ttilde{K})$ holds; (see curve
$\gamma_2$ in Figure~\ref{fig:sigma0}); \\
{\it (iii)} if $K_1 \in [0, \widetilde{K}_1)$, then $z(t) > 0$, before
the trajectory intercept $\Sigma_0$ and $z(t) < 0$ after that (see curve
$\gamma_3$ in Figure~\ref{fig:sigma0}).

\begin{folge} \label{kor:nullsig}                                   
Let $\widetilde{K}_1$, $\ttilde{K}$, $\breve{x}$ and $l$ be defined as in
Lemma~\ref{lem:nullst} and in the proof of Corollary~\ref{kor:nullkurve}.
Then there exists $\hat{x} \in (\breve{x},\tilde{x})$, $\widehat{K} > 0$
and a continuous mapping $h_0: [\hat{x}, \tilde{x}] \to [\widehat{K},
\ttilde{K}]$ with
$$  \Sigma_0 \cap \{ (x_0,K_0)\, |\, F(x_0) \le K_0 x_0 \} =
    \{ (x_0,h_0(x_0))\, |\, x_0 \in [\hat{x},\tilde{x}] \} \, . $$
Moreover:
\begin{itemize}
\item [(a)] $h_0(\hat{x}) = \widehat{K}$, $h_0(\tilde{x}) = \ttilde{K}$,
            \ $F(\hat{x}) = \widehat{K} \hat{x}$;
\item [(b)] $h_0$ is continuous differentiable in $(\hat{x},\tilde{x})$
            and $h_0'(x) > 0$, $x \in [\hat{x},\tilde{x})$.
\end{itemize}
\bew \rm \
Consider the mapping
$$ \xi: [0,\widetilde{K}_1] \ni K_1 \auf ( x(l(K_1);K_1), z(l(K_1);K_1) )
   \in [0,\tilde{x}] \times \R \, , $$
where $x(\cdot;\cdot)$ and $z(\cdot;\cdot)$ are defined as in the proof of
Corollary~\ref{kor:nullkurve}. Then
$$ \xi'(K_1) = 
\Big(  x'(l(K_1);K_1)l'(K_1) +
       \frac{\partial x(\cdot;\cdot)}{\partial K_1}(l(K_1);K_1) , 
       z'(l(K_1);K_1)l'(K_1) +
       \frac{\partial z(\cdot;\cdot)}{\partial K_1}(l(K_1),K_1) \Big) $$
and we see by the implicit function theorem that the region $\Sigma_0
\cap \{ (x_0,K_0)\, |\, F(x_0) \le K_0x_0 \}$ may be reparameterized by
a function $h_0$. Notice that the implicit function theorem may be used
due to Corollary~\ref{kor:nullkurve}. The condition $h_0'(x) > 0$ follows
from the same corollary.
\eop
\end{folge}

\begin{figure}[!ht]
\begin{center}
\epsfig{figure=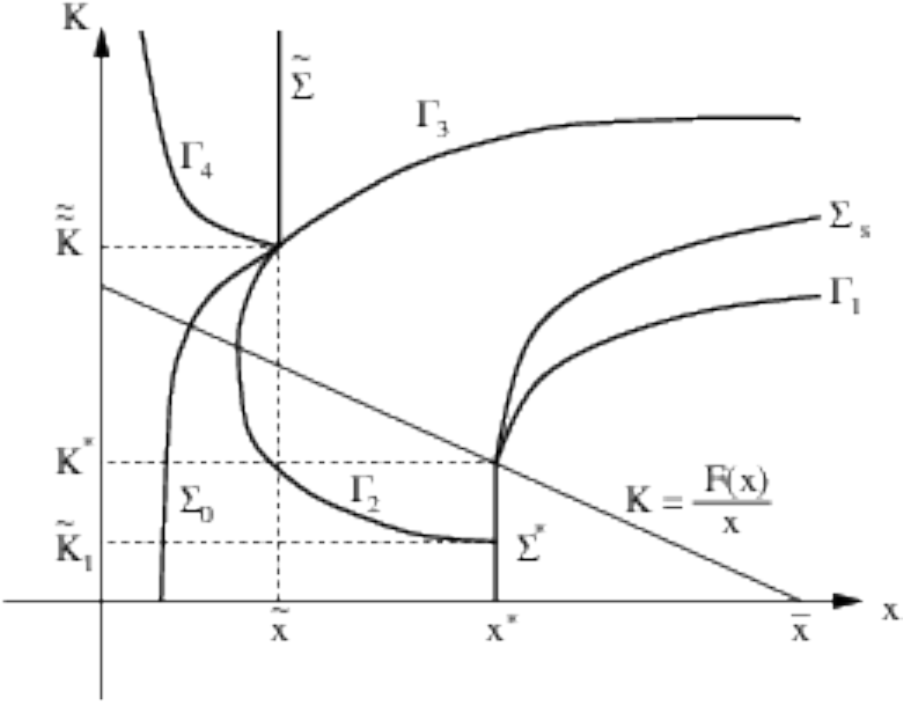,width=.75\textwidth}
\caption{Main curves.} \label{fig:szenario}
\end{center}
\end{figure}

Let $(x,K)$ be the solution of the initial value problem
\begin{equation} \label{eq:gam3}
\left\{ \begin{array}{l@{}l}
   x' = -F(x) + K x & ,\ x(0) = \tilde{x} , \\
   K' = \gamma K    & ,\ K(0) = \ttilde{K}
\end{array} \right.
\end{equation}
in the interval $[0,\tau)$. Similar to the definition of $\Gamma_1$ we
obtain $\bar{t} \in (0,\tau)$ with $x(\bar{t}) = \bar{x}$. We denote
the curve
$$ [0,\bar{t}] \ni t \longmapsto
   (x(t),K(t)) \in [\tilde{x},\bar{x}] \times [\ttilde{K},\infty) $$
by $\Gamma_3$. (In Figure~\ref{fig:sigma0} one can recognize $\Gamma_3$ as
the part of the curve $\gamma_2$ with $x > \tilde{x}$ and $K > \ttilde{K}$.)
Now let $(x,K)$ be the solution of the initial value problem
\begin{equation} \label{eq:gam5}
\left\{ \begin{array}{l@{}l}
   x' = F(x) - K x & ,\ x(0) = \tilde{x} \, , \\
   K' = -\gamma K  & ,\ K(0) = \ttilde{K} \, .
\end{array} \right.
\end{equation}
This solution meets the curve $\Sigma^*$ in $(x^*,\widetilde{K}_1)$ for
some $\tau>0$ (see Lemma~\ref{lem:nullst}). The curve defined by this
trajectory is called $\Gamma_2$. (In Figure~\ref{fig:sigma0}, the curve
$\Gamma_2$ corresponds to the part of $\gamma_2$ with $K < \ttilde{K}$.)
The last curve we define is $\Gamma_4$, which is parameterized by the
solution $(x,K)$ of
\begin{equation}
\left\{ \begin{array}{l@{}l}
   x' = -F(x)    & ,\ x(0) = \tilde{x} \, , \\
   K' = \gamma K & ,\ K(0) = \ttilde{K} \, .
\end{array} \right.
\end{equation}
Note that the solution exists in $[0,\infty)$ and $x'(0) < 0$, $K'(0) > 0$
(see Figure~\ref{fig:szenario}).

Now we come back to the region $x > x^*$, more specifically, above the
curve $\Gamma_1$. We already know that if the initial condition $(x_0,K_0)$
lays below $\Gamma_1$, then $K_{0,+} \ge h_1(x_0)$ must hold. We want to
determine a curve $\Sigma_s$, above $\Gamma_1$, upon which the trajectories
must jump for this initial conditions, i.e. $(x_0,K_{0,+}) \in \Sigma_s$.
In order to be able to construct such a curve $\Sigma_s$, we need the
fact that there exists $(x_0,K_0) \in (x^*,\bar{x}) \times (K^*,\infty)$
such that the solution $(x,K,z,\la)$ of
\begin{eqnarray*}
\left\{ \begin{array}{l@{}l}
  x' = F(x) - K x                & ,\ x(0) = x_0 \, , \\
  K' = -\gamma K                 & ,\ K(0) = K_{0,+} \, , \\
  z'   = (z-r') g(x) - \psi_*(x) & ,\ z(0) = z_0 \, , \\
  \la' = (\la -r)\kappa - z + r' & ,\ \la(0) = r \, ,
\end{array} \right.
\end{eqnarray*}
meets the curve $\Sigma^*$.

\begin{lemma} \label{lem:zurueck}                                   
There exists $\widetilde{K}_2 \in (0,\widetilde{K}_1)$, $\widetilde{K}_3
\in (\widetilde{K}_1,K^*)$ and $a > 0$ such that the solution
$( x(\cdot;K_1), z(\cdot;K_1), \la(\cdot;K_1) )$ of the system
\begin{equation} \label{eq:sys_xzla}
\left\{ \begin{array}{l@{}l}
  x'   = -F(x) + K_1e^{\gamma t} x & ,\   x(0) = x^* \, , \\
  z'   = -(z-r') g(x) + \psi_*(x)  & ,\   z(0) = r' \, , \\
  \la' = - (\la -r)\kappa + z - r' & ,\ \la(0) = r \, ,
\end{array} \right.
\end{equation}
exists in $[0,a]$ for each $K_1 \in (\widetilde{K}_2,K^ *)$. Moreover, for
each $K_1 \in (\widetilde{K}_3,K^ *)$ there exist numbers $\rho(K_1)$,
$\sigma(K_1)$, $\tau(K_1)$ with
\begin{itemize}
\item [(a)] $0 < \rho(K_1)<\sigma(K_1)<\tau(K_1)< a$;
\item [(b)] $x(t;K_1) \le x^ *$, $t \in (0,\rho(K_1))$, $x(\rho(K_1);K_1) =
            x^*$, $x(t;K_1) > x^*$, $t \in (\rho(K_1),a]$, \\
            and $x(a;K_1) \ge \bar{x}$;
\item [(c)] $0 < z(t;K_1)<r'$, $t \in (0,\sigma(K_1))$, $z(\sigma(K_1);K_1)
            = r'$, $z(t;K_1) > r'$, $t \in (\sigma(K_1),a]$;
\item [(d)] $\la(t;K_1) < r$, $t \in (0,\tau(K_1))$, $\la(\tau(K_1);K_1) = r$,
            $\la'(\tau(K_1);K_1) > 0$;
\item [(e)] $\D\lim_{K_1 \uparrow K^*} \tau(K_1) = 0$.
\end{itemize}
\bew \rm \
Since the trajectory for $K_1 := K^*$ corresponds to $\Gamma_1$ we can
choose $a > 0$ and $\widetilde{K}_2 \in (0,\widetilde{K}_1)$ with
$x(a;K_1) \ge \bar{x}$, $K_1 \in (\widetilde{K}_2,K^*)$. It follows
from the differential equations for $z$ and $\la$ that $z(t;K^*) > r'$,
$\la(t;K^ *) > r$, $t \in (0,a]$. \\
Let $\eps \in (0,a)$ be given. Choose $\beta > 0$ and $\widetilde{K}_3 \in
(\widetilde{K}_1,K^*)$ such that for each $K_1 \in (\widetilde{K}_3,K^*)$
$$ z(t;K_1) \ge r' + \beta \, ,\ \la(t;K_1) \ge r + \beta \, ,\
   t \in [\eps,a] \, . $$
This can be done due to the fact that the solution of \eqref{eq:sys_xzla}
depends continuously on the parameter $K_1$. \\
Let $K_1 \in (\widetilde{K}_3,K^*)$ and set $(x,z,\la) := (x(\cdot;K_1),
z(\cdot;K_1), \la(\cdot;K_1))$. Then we see that there exists $\rho(K_1) > 0$
with the property in {\it (b)}. \\
Since $z''(0) = \psi'_*(x^*)\, x'(0) < 0$ and since $z(a) \ge r' + \beta$,
we obtain $\sigma(K_1) > 0$ with $z(\sigma(K_1)) = r'$ and $0 < z(t) < r'$,
$t \in (0,\sigma(K_1))$. $\sigma(K_1) < \rho(K_1) $ cannot hold since
$z'(\sigma(K_1)) = \psi_*(x(\sigma(K_1))) < 0$. $\sigma(K_1) = \rho(K_1)$
cannot hold since $z''(\sigma(K_1)) = z''(\rho(K_1)) = \psi_*(x^*)
x'(\rho(K_1)) > 0$. From the differential equation for $z$ and the fact
that $\rho(K_1) < \sigma(K_1)$ we conclude $z(t) > r'$, $t > \sigma(K_1)$.
Repeating the argumentation above we obtain $\tau(K_1)$ with
$\la(\tau(K_1)) = r$, $\la(t) < r$, $0 < t < \tau(K_1)$. $\tau(K_1) <
\sigma(K_1)$ cannot hold since $\la'(\tau(K_1)) = z(\tau(K_1)) - r' < 0$.
Assume $ \tau(K_1) = \sigma(K_1)$. Then $\la'(\tau(K_1)) = z'(\sigma(K_1))
= \psi_*(x(\sigma(K_1))) > 0$. This contradicts the fact that $\la(t) < r$,
$t \in (0,\tau(K_1))$. From the differential equation for $\la$ and the
fact that $z(t) > r'$, $t \in \tau(K_1)$, we obtain $\la(t) > r$,
$t > \tau(K_1)$. Clearly, from the continuous dependency we obtain
$\D\lim_{K_1 \uparrow K^*} \tau(K_1) = 0$.
\eop
\end{lemma}

The value $\tau(K_1)$ according to Lemma~\ref{lem:zurueck} is locally
uniquely determined and the same is true for $K(\tau(K_1))$. We want to
identify $K(\tau(K_1))$ as a value $K_{0,+}$, when an extremal trajectory
starts in the initial value  $(x_0,K_0) \in (x^*,\bar{x}) \times [0,\infty)$,
with $K_0 < h_1(K_0)$. To do this we need more information concerning the
mapping $K_1 \mapsto \tau(K_1)$.

Consider the system \eqref{eq:sys_xzla} for $K_1 \in (\widetilde{K}_3,K^*)$
and let $(x,\la,z)$ be the corresponding solution. To make clear the
dependence of $\la$ on $K_1$, we denote it by $G(\cdot;K_1)$. The equation
\begin{equation} \label{eq:glei} 
G(\tau;K_1) = r
\end{equation}
describes the fact that the solution $G(\cdot;K_1)$ has the value $r$ at
time $\tau$. In order to find the jump curve ${\it \Sigma}_s$ we try to
resolve (\ref{eq:glei}) with respect to $\tau$. Note that from Lemma~%
\ref{lem:zurueck}, we know that there are solutions of \eqref{eq:glei}.

\begin{lemma} \label{lem:regu}                                      
Using the notations of Lemma~\ref{lem:zurueck} the following assertions hold:
\begin{itemize}
\item [(a)] There exists $\widehat{K}_1 \in [\widetilde{K}_1,K^ *)$
            and a continuous differentiable mapping $g_{\tau}:
            (\widehat{K}_1,K^*) \to (0,\infty)$ with
            \begin{equation} \label{eq:glei1} 
            G(g_{\tau}(K_1);K_1) = r,\ K_1 \in (\widehat{K}_1,K^*) ;
            \end{equation}
\item [(b)] $g_{\tau}'(K_1) < 0$ for all $K_1\in (\widehat{K}_1,K^*)$;
\item [(c)] $x(g_{\tau}(\widehat{K}_1);\widehat{K}_1) \in \Gamma_3$ or
            $\widehat{K}_1 = \widetilde{K}_1$ and
            $x(g_{\tau}(\widehat{K}_1);\widehat{K}_1) = \bar{x}$.
\end{itemize}
\bew \rm \
Let $K_1 \in (\widetilde{K_3}, K^ *)$. Since $\la = G(\cdot;K_1)$, the
function $G$ is obviously differentiable and we have
$$ \frac{\partial G}{\partial \tau}(\tau(K_1);K_1) = \la'(\tau(K_1)) > 0 $$
(see Lemma~\ref{lem:zurueck}). With the implicit function theorem we obtain
a neighborhood $U$ of $K_1$ and a function $g_{\tau}: U \to (0,\infty)$
such that $G(g_{\tau}(K_1);K_1) =r$, $K_1 \in U$ holds. The implicit
function theorem implies more: $g_{\tau}$ can be extended in a maximal
way to $(\widehat{K}_1,K^*)$ with $\widehat{K}_1 \in [\widetilde{K}_1,K^*)$
and \eqref{eq:glei1} holds. The properties in {\it c)} are a consequence of
the maximality of the extension. Moreover:
$$ \la'(g_{\tau}(K_1)) g_{\tau}'(K_1) = - v(g_{\tau}(K_1)) \, , $$ 
where $v:= \frac{\D \partial \la}{\D \partial K_1}(\cdot)$ solves
$$ v' = -\kappa v + w \, ,\ v(0)= 0 \, , $$
$w:= \frac{\D \partial z}{\D \partial K_1}(\cdot)$ solves
$$ w'= -w g(x(t)) + \big[-(z(t)-r')g'(x(t)) + \psi_*'(x(t))\big] y \, ,\
   w(0) = 0 \, , $$
and $y := \frac{\D \partial x}{\D \partial K_1}(\cdot)$ solves
$$ y' = \big[-F'(x(t)) + K_1 e^{\gamma t}\big] y + e^{\gamma t} x(t) \, ,\
   y(0) = 0 . $$
Obviously $y(t) > 0$, $t \in (0,g_{\tau}(K_1)]$. Since $-(z(t)-r') g'(x(t))
+ \psi_*'(x(t)) > 0$ for each $t \in (0,g_{\tau}(K_1)]$, we have $w(t) > 0$,
$t \in (0,g_{\tau}(K_1)]$, and therefore $v(t) = \frac{\D \partial \la}
{\D \partial K_1}(t) > 0$, $t \in (0,g_{\tau}(K_1)]$. It follows
$g_{\tau}'(K_1) < 0$.
\eop
\end{lemma}

Now, we have to distinguish two cases: \\[1ex]
\begin{tabular}{ll}
\textbf{Case I} & $x(g_{\tau}(\widehat{K}_1);\widehat{K}_1) \in \Gamma_3$ and
                  $x(g_{\tau}(\widehat{K}_1);\widehat{K}_1) < \bar{x}$; \\[1ex]
\textbf{Case II}& $x(g_{\tau}(\widehat{K}_1);\widehat{K}_1) = \bar{x}$.
\end{tabular} \\[1ex]
In each case we have a curve $\Sigma_s$ defined by 
$$ [\widehat{K}_1,K^ *] \ni K_1 \auf ( x(g_{\tau}(K_1);K_1), K_1 e^{\gamma
   g_{\tau}(K_1)} ) \in [x^*,\bar{x}] \times [K^*,\infty) . $$
In case II this curve ends at the \glqq boundary\grqq $x = \bar{x}$. In
case I we want to continue this curve such that the continuation ends
at the \glqq boundary\grqq $x = \bar{x}$ too. For this continuation
the trajectories starting in $(\tilde{x},K)$, $K \ge \ttilde{K}$, come
into consideration. To analyze the situation we need the curve
$\widetilde{\Sigma}$ which is defined by
$$ \widetilde{\Sigma}: [0,\infty) \ni K \auf (\tilde{x},K) \in
   [0,\bar{x}] \times [0,\infty)\, . $$
(See Figure~\ref{fig:szenario}). As we will see in Theorem \ref{th:sigma0}
an optimal trajectory $(x,K)$ with adjoint variables $(z,\la)$ meets the 
curve $\widetilde{\Sigma}$ with the values
$$ z(\sigma)=0\, ,\ z'(\sigma)=0\, ,\ \la(\sigma) < r . $$
This motivates the construction of the continuation of $\Sigma_s$ in the
following way: Compute trajectories $(x,K,z)$ backwards in time starting
with initial values
$$ x(0) = \tilde{x}\, ,\ K(0) = K_0 > \ttilde{K}\, ,\ z(0) = 0 . $$

\begin{lemma}\label{lem:case2}                                      
There exists $\ttilde{K_1} > \ttilde{K}$ and $a > 0$ such that for each
$K_1 \ge \ttilde{K_1}$
$$ x(a;K_1) = \bar{x}\, ,\ z(t;K_1)< r' \,,\ t \in [0,a] , $$
where $(x(\cdot;K_1),z(\cdot;K_1))$ is the solution of 
\begin{equation}\label{eq:sys_xzc1}
\left\{ \begin{array}{l@{}l}
  x'   = -F(x) + K_1e^{\gamma t} x & ,\   x(0) = \tilde{x} , \\
  z'   = -(z-r') g(x) + \psi_*(x)  & ,\   z(0) = 0\,.  
\end{array} \right.
\end{equation}
\bew \rm \
Let $\eps > 0$. It follows from the differential equation for $x$ that the
solution $x(\cdot;K_1)$ reaches $x = \bar{x}$ for a time $t_1 < \eps$ if
$K_1$ is sufficiently large. Since the differential equation for $z$ may
be considered as a linear equation (if we plug in $x(\cdot;K_1)$) the
value $z=r'$ cannot be reached in the time interval $[0,\eps]$ if $K_1$
is sufficiently large.
\eop
\end{lemma}

As we know from the results above an optimal trajectory $(x,K)$ with adjoint
variables $(z,\la)$ starts in $(\tilde{x},\ttilde{K})$ in the following way:
$$ x(0) = \tilde{x}\, ,\ K(0) =\ttilde{K}\, ,\ z(0)=0\, ,\ z'(0) = 0\, , \
   \tilde{\la} := \la(0) < r . $$
For each $K_1 > \ttilde{K}$ there exists a time $\xi =\xi(K_1)$ with
$K_1 e^{-\gamma \xi} = \ttilde{K}$. Set $\Lambda(K_1) := \la(\xi;K_1)$
where $\la(\cdot;K_1)$ is the solution of 
$$  \la' = -(\la-r) \kappa - r' \, ,\ \ \la(0) = \tilde{\la} . $$
Now consider for each $K_1 > \ttilde{K}$ the solution of
\begin{equation}\label{eq:sys_xzlc1}
\left\{ \begin{array}{l@{}l}
  x'   = -F(x) + K_1e^{\gamma t} x & ,\   x(0) = \tilde{x}, \\
  z'   = -(z-r') g(x) + \psi_*(x)  & ,\   z(0) = 0 , \\
  \la' = - (\la -r)\kappa + z - r' & ,\ \la(0) = \Lambda(K_1)\, .
\end{array} \right.
\end{equation}
We want to find for $K_1 > \ttilde{K}$ some time $\tau(K_1)$ such that
$\la(\tau(K_1);K_1)= r$ holds. Again we use the notation $G(\cdot;K_1)
:= \la(\cdot;K_1)$.

\begin{lemma}\label{lem:gc1}                                       
In case I the following assertions hold:
\begin{itemize}
\item [(a)] There exists $\ttilde{K_1} \in (\ttilde{K},\infty)$ 
and a continuous differentiable mapping $g_{\tau}: (\ttilde{K},\ttilde{K_1}) \to(0,\infty)$ with
            \begin{equation}\label{eq:glei2} 
            G(g_{\tau}(K_1);K_1) = r, K_1 \in (\ttilde{K},\ttilde{K_1}) .
            \end{equation}
\item [(b)] $g_{\tau}'(K_1) < 0$ for all $K_1\in (\ttilde{K},\ttilde{K_1})\,.$
\item [(c)] $x(g_{\tau}(\ttilde{K_1});\ttilde{K_1})  = \bar{x}\,.$
\end{itemize}
\bew \rm \
This can be proved similar to Lemma \ref{lem:regu}. The main observation
is that we have $z(\tau;K_1)>r'$ if $\la(\tau;K_1)=r$. The result of
{\it (c)} is a consequence of Lemma~\ref{lem:case2}.
\eop
\end{lemma}

Now we have constructed a curve $\Sigma_s$ connecting $(x^*, K^*)$ with
$(\bar{x}, \bar{K})$, where $\bar{K} := \ttilde{K_1} e^{\gamma g_{\tau}
(\ttilde{K_1})}$.
Since this curve should be used as a jump curve in the region
$[x^*,\bar{x}] \times [0,\infty)$, we want to reparameterize this curve in
such a way that $[x^*,\bar{x}]$ is the parameter interval. But to do this
we need the fact that the following function
$$ [\widehat{K}_1,K^ *] \cup [\ttilde{K},\ttilde{K_1}] \ni K_1 \auf
   x(g_{\tau}(K_1);K_1) \in [x^ *,\bar{x}] $$
is monotone increasing. Unfortunately, we are not able to prove this.
The fact we can prove is that the mapping
$$ H : [\widehat{K}_1, K^*] \ni K_1 \auf x(g_{\tau} (K_1); K_1) \in
   [x^*, \bar{x}] $$
is monotone increasing in a neighborhood of $K^*$.

\begin{lemma}\label{lem:taup}                                       
There exists $K_+ \in [\widehat{K}_1, K^*)$ and $m_0,\ m_1 > 0$ such that
\begin{itemize}
\item[a)] $\tau (K^*) = 0$;
\item[b)] $\tau'(K^*) = -\frac{\D 4}{\D \gamma K^*}$, $\tau$ is
          differentiable in $[K_+, K^* ]$;
\item[c)] $-m_1 (K^* - K) \le \tau(K) \le -m_0 (K^* - K)$ for all
          $K \in [K_+, K^*]$.
\end{itemize}
\bew \rm \
Item {\it a)} is obvious. Clearly, $ \tau $ is differentiable in $[K_+,K^*)$.
From the identity
$$ \la (\tau(K); K) = 0, K \in [\widehat{K}_1, K^*], $$
we obtain
$$ \la' (\tau(K); K) \tau'(K) + v(\tau(K); K) = 0,\
   K \in [\widehat{K}_1, K^*). $$
Here we have used the notation of Lemma \ref{lem:zurueck}. Due to
the fact that $\la'(\tau(K); K) > 0$ for each $K \in [\hat{K}_1,K^*)$
we have
$$ \tau' (K) =- \frac{\D  v(\tau(K); K)}{\D \la'(\tau (K); K)} ,\
   K \in [\widehat{K}_1, K^*). $$
Using Taylor's expansion for $v$ and $\la$ we obtain
$$ \tau'(K) = \frac{\D -\frac{\D 1}{\D 6} v^{(''')} ( \xi; K)}
   {\D \frac{\D 1}{\D 24} \la^{(iv)} ( \eta; K)} ,\
   K \in [\widehat{K}_1,K^*) , $$
where $ \xi, \eta, \in (0, \tau(K))$. From this we can conclude
$$ \tau'(K^*) =-4 \frac{\D v^{(''')}(0; K^*)}{\D \la^{(iv)} (0; K^*)}
   = -\frac{\D 4}{\D \gamma K^*} . $$
The result in {\it c)} is a consequence of the estimate in {\it b)} by
using continuity arguments.
\eop
\end{lemma}

\begin{lemma}\label{lem:hmon}                                       
There exists $ K_+ \in [\hat{K}_1, K^*) $ such that
$$ H(K^*) = x^*, H'(K^*) = 0; H' (K) > 0, K \in (K_+, K^*). $$
\bew \rm \
Let $ K_+$ be chosen as in Lemma \ref{lem:taup}.
We have 
$$ H'(K) = x' (\tau(K); K) \tau'(K) + y (\tau(K); K), K \in [K_+, K^*], $$
and by Lemma \ref{lem:taup} we have $ H' (K^*) = 0$. From the differential
equation for $y$ we obtain that there exists $M > 0$ such that
$$ y (\tau (K); K) \ge M \tau (K) , K \in [K_+, K^*].$$
Since $\tau'$ is bounded in $[K_+, K^*]$ and since $x'(\tau(K); K) = 0$
we obtain the assertion, eventually by making $K_+$ larger.
\eop
\end{lemma}

\begin{folge} \label{kor:sigmah}                                    
There exists $x_s \in (x^*, \bar{x}]$ and a continuously differentiable
mapping $h_s :[x^*, x_s]$ such that
$$ (x,h_s(x)) \in \Sigma_s\, ,\ x \in [x^*, x_s] . $$
\bew \rm \
The mapping $h_s$ can be found by reparameterizing the curve
$$ K_1 \auf (x(g_{\tau} (K); K) , K(g_{\tau} (K); K)) $$
using the results of Lemma \ref{lem:hmon}.
\eop
\end{folge}

Now the jump curve $\Sigma_s$ allow us to find, given $(x_0,K_0)$ with
$x_0 \in (x^*,x_s)$ and $K_0 < h_s(x_0)$, the optimal initial value
$K_{0,+} := h_s(K_0)$. Notice that in the region around $\Gamma_1$,
$\Sigma^*$, and $\Sigma_s$, the behavior of the extremals is rather
clear. Notice also that $h_s(x) > h_1(x)$, for $x \in (x^*,\bar{x})$,
since $\Sigma_s$ is above $\Gamma_1$. From the construction of $\Sigma_s$
and $\Gamma_3$, we may have two cases:
$\Sigma_s$ and $\Gamma_3$ have a point in common; $\Sigma_s$ and $\Gamma_3$
do not intersect. In the next section we are able to find for each initial
condition $(x_0,K_0) \in (0,\bar{x}) \times (0,\infty)$ the corresponding
optimal trajectory for the case where $\Sigma_s$ and $\Gamma_3$ do not
intersect; the other case can be handled in an analogous way.
%
%
%
\section{Optimal trajectories}\label{sec:optim}

\begin{figure}[!ht]
\begin{center}
\epsfig{figure=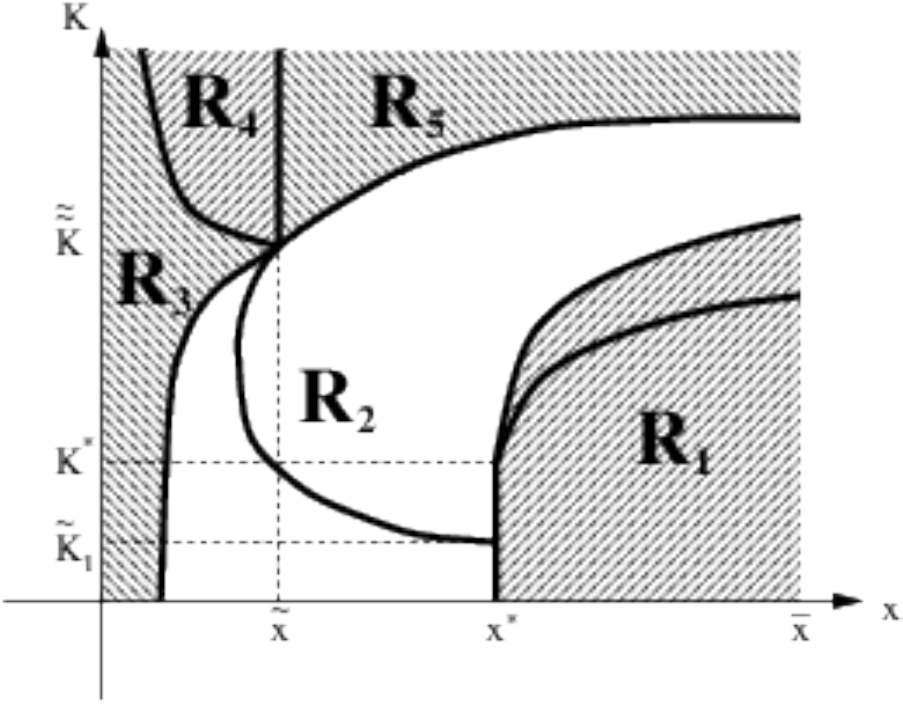,width=.75\textwidth}
\caption{Main regions.} \label{fig:regions}
\end{center}
\end{figure}

In this section we summarize the previous results, in order to design a 
complete picture of the optimal trajectories. We are also able to determine 
the correspondent optimal controls, based on the information we obtain from 
the switch variables. Let us consider the regions defined by
$$  \Sigma^*, \widetilde{\Sigma}, \Sigma_s, \Sigma_0, \Gamma_1, \Gamma_2,
    \Gamma_3, \Gamma_4 $$
in $[0,\bar{x}] \times [0,\infty)$. In Figure \ref{fig:regions} we present a 
sketch of the five main regions for the case of the logistic
function. Figure \ref{fig:regions} shows the case that $ x_s \ge \bar{x} $
holds. The analysis in the following is done under the following assumption:
$$ x_s = \bar{x}\, . $$
Unfortunately, we are not able to present reasonable conditions which
imply this assumption.

\begin{description}
\item [Domain (R1):] boundaries $\Sigma^*$, $\Sigma_s$, $\{(x,K) \in 
[0,\bar{x}] \times [0,\infty)\ |\ x= \bar{x}\}$;
\item [Domain (R2):] boundaries $\Sigma_0$, $\Sigma^*$, $\Sigma_s$, $\Gamma_3$ 
and eventually $\{(x,k) \in [0,\bar{x}] \times [0,\infty)\ |\ x = \bar{x}\}$;
\item [Domain (R3):] boundaries $\{(x,K) \in [0,\bar{x}] \times [0,\infty)\ 
|\ x= 0\}$, $\Sigma_0$, $\Gamma_4$;
\item [Domain (R4):] boundaries $\{(x,K) \in [0,\bar{x}] \times [0,\infty)\ 
|\ x= 0\}$, $\Gamma_4$, $\widetilde{\Sigma}$;
\item [Domain (R5):] boundaries $\widetilde{\Sigma}$, $\Gamma_3$, $\{(x,K) 
\in [0,\bar{x}] \times [0,\infty)\ |\ x = \bar{x}\}$.
\end{description}

\begin{rem}\label{rem:erl}                                            
In the following we may assume, without loss of generality, that $K_0 = 
K_{0,+}$ holds since in the case of $K_0 < K_{0,+}$ we may leave the region 
where we started or not. In the first case we have to apply the discussion 
in another region, in the second case we have to repeat the analysis in the 
same region with $K_0 := K_{0,+}$.
\eor
\end{rem}

\begin{satz}\label{th:rrr}                                           
Let $(x,K,u,\mu)$ be optimal with adjoint variables $z$, $\la$. The following 
assertions hold:
\begin{itemize}
\item [(a)] Let $(x_0,K_0)$ be in (R1). Then $z_0 > r'$, $\la_0 = r$, 
$(x_0,K_{0,+}) \in \Sigma_s$ and there exists $\tau > 0$ with $(x(\tau), 
K(\tau)) \in \Sigma^*$ and $u \equiv 1$, $\mu \equiv 0$ in $(0,\tau)$;
\item [(b)] Let $(x_0,K_0)$ be in (R2). Then $z_0 > 0$, $\la_0 < r$ and there 
exists $\tau > 0$ with $(x(\tau),K(\tau)) \in \Sigma^*$ and $u \equiv 1$, 
$\mu \equiv 0$ in $(0,\tau)$;
\item [(c)] Let $(x_0,K_0)$ be in (R3). Then $z_0 < 0$, $\la_0 < r$ and there 
exists $\tau > 0$ with $(x(\tau),K(\tau)) \in \Sigma_0$ and $u \equiv 0$, 
$\mu \equiv 0$ in $(0,\tau)$;
\item [(d)] Let $(x_0,K_0)$ be in (R4). Then $z_0 < 0$, $\la_0 < r$ and there 
exists $\tau > 0$ with $(x(\tau),K(\tau)) \in \widetilde{\Sigma}$ and 
$u \equiv 0$, $\mu \equiv 0$ in $(0,\tau)$;
\item [(e)] Let $(x_0,K_0)$ be in (R5). Then $z_0 > 0$, $\la_0 < r$ and there 
exists $\tau > 0$ with $(x(\tau),K(\tau)) \in \widetilde{\Sigma}$ and 
$u \equiv 1$, $\mu \equiv 0$ in $(0,\tau)$.
\end{itemize}
\bew \rm \
Ad {\it (a)}: We actually prove only that if $(x_0,K_0)$ is in $(R1)$, then 
$K_{0,+} \ge h_s(x_0)$, i.e. we must jump either to (R2) or to (R5). Later 
on, in the proof of items (b) and (e), we will see that the initial condition 
$\la_0 = r$ is not allowed in these regions. The last possible case: 
$(x_0,K_{0,+}) \in \Gamma_3$ is excluded in Theorem~\ref{th:gamma} \\
From Lemma~\ref{lem:sprung3}, it is enough to consider initial conditions 
$(x_0,K_0)$ with $K_0 \ge h_1(x_0)$. Assume $\la_0 < r$. Then there exists 
a $\tau > 0$ with $\la(\tau) = r$ and $\la(t) < r$, $t \in [0,\tau)$.
Consequently $\la'(\tau)=0$, $z(\tau)=r'$, $z'(\tau) = -\psi_*(x(\tau))$. \\
We consider tree cases: {\it i)} $z'(\tau) < 0$: then $\la''(\tau) > 0$,
contradicting (R).
{\it ii)} $z'(\tau) > 0$: then $x(\tau) < x^*$ and $(x(\tau), K(\tau)) \in
\cup_{i=2}^5 (Ri) \cup \widetilde{\Sigma} \cup \Sigma_0 \cup \Gamma_2 \cup
\Gamma_3 \cup \Gamma_4$. As we will see, the initial condition $\la_0 = r$
is not allowed in these regions (curves), and again we have a contradiction.
{\it iii)} $z'(\tau) = 0$: then $x(\tau) = x^*$. From Lemma~\ref{lem:lrplus},
$K(\tau) > K^*$ cannot occur and we must have $(x(\tau), K(\tau)) \in
\Sigma^*$. From the differential equation for $(x,K)$, follows the existence
of $\sigma \in (0,\tau)$ such that $(x(\sigma), K(\sigma)) \in \Sigma_s$.
However, $\la(\sigma) < r$, contradicting the construction of $\Sigma_s$,
since the optimal trajectory $(x,K)$ hits the curve $\Sigma^*$. \\
Therefore, we must have $\la_0 = r$, what implies $z_0 \ge r'$. Note
that $z_0 = r'$ is not possible, since we would have $\la'(0) = 0$ and
$\la''(0) = -z'(0) = \psi_*(x_0) > 0$, contradicting (R). Finally, we
exclude the case $K_{0,+} < h_s(x_0)$. If this where not the case, we
would obtain a contradiction arguing as in the case $\la_0 < r$ above.
\medskip

\noindent
Ad {\it (b)}: Assume $\la(0) = r$. Then $\la'(0) \le 0$ and $z_0 \ge r'$. \\
If $z_0 = r'$, then $\la'(0) = 0$ and we have three possible cases:
{\it i)} $x_0 > x^*$: we have $\la''(0) = \psi_*(x_0) > 0$, contradicting (R);
{\it ii)} $x_0 = x^*$: cannot occur, due to Lemma~\ref{lem:lrplus};
{\it iii)} $x_0 < x^*$: then $z'(0) > 0$ and $\la''(0) = -z'(0) < 0$. Since
$\la(t) < r$ for all $t > 0$ is not possible, there exists $\tau > 0$ with
$\la(\tau) = r$, $\la(t) < r$, $t \in (0,\tau)$. Then $\la'(\tau) = 0$,
$z(\tau) = r'$ and there exists $\sigma \in (0,\tau)$ with $z'(\sigma) = 0$,
$z(\sigma) > r'$. From the differential equation for $z$ we conclude with
(V4) that $x(\tau) \ge x(\sigma) > x^*$. This is a contradiction to
Lemma~\ref{lem:sprung1}, since it is easy to see that $x(\tau) K(\tau)
< F(x(\tau))$.

\noindent
For $z_0 > r'$ we have again three possible cases:
{\it i)} $x_0 > x^*$; {\it ii)} $x_0 = x^*$; {\it iii)} $x_0 < x^*$.
Cases {\it ii)} and {\it iii)} are excluded analogous as above.
In case {\it i)}, since $\la'(0) < 0$, there exists $\tau > 0$ such
that $\la(\tau) = r$, $\la'(\tau) = 0$, $z(\tau) = r'$, $z'(\tau) =
-\la''(\tau) \ge 0$. Therefore, we have $x(\tau) \le x^*$. If $x(\tau)
< x^*$, follows from $(x_0, K_0) \in (R_2)$ and the fact that there
are no jumps in interval $(0,\tau)$ that $(x(\tau), K(\tau)) \in (R_2)$
must hold. In this case a contradiction can be obtained arguing as in
the case $\la(0) = r$ \& $z(0) = r'$ above.
If $x(\tau) = x^*$, follows from Lemma~\ref{lem:lrplus} that $K(\tau) >
K^*$ cannot occur. Thus, we must have $(x(\tau), K(\tau)) \in \Sigma^*$.
However, from the construction of $\Sigma_s$, we know that along every
optimal arc that hits $\Sigma^*$ (starting at $(x_0,K_0) \in$ (R2) with
$x_0 \ge x^*$) we must have $\la(t) < r$, $t \in [0,\tau)$. In particular,
$\la(0)<r$ must hold, which is a contradiction.
Therefore, we have $\la_0 < r$. 

Since $\la(t) < r$ for all $t > 0$ is not
allowed, there exists $\tau > 0$ with $\la(\tau) = r$, $\la(t) < r$,
$t \in [0,\tau)$, $\la'(\tau) = 0$, $z(\tau) = r'$.
Since the initial condition $\la(0) = r$, is not allowed in (R2),
we conclude that the trajectory must leave (R2) at some time $\sigma \le
\tau$. Since there are no jumps in the interval $[0,\tau)$, the
trajectory can leave (R2) only through $\Sigma_s$ or $\Sigma^*$. If
$(x(\sigma), K(\sigma)) \in \Sigma_s$ we have two possibilities:
{\it i)} $\sigma < \tau$: in this case we can find $\eps > 0$ such that
$(x(\sigma+\eps), K(\sigma+\eps)) \in$ (R1) and $\la(\sigma+\eps) < r$.
From item {\it (a)} above we know that this cannot occur. {\it ii)}
$\sigma = \tau$: in this case we have $(x(\tau), K(\tau)) \in
\Sigma_s$. However, this is not in agreement with the inequality
$x(\tau) \le x^*$, which follows from $z'(\tau) = -\la''(\tau) \ge 0$
and $z(\tau) = r'$. Therefore, the trajectory must leave (R2) through
$\Sigma^*$. If $\sigma < \tau$, then $\la(\sigma) < r$ and we have a
contradiction by Lemma~\ref{lem:sprung3}.
Thus the trajectory must leave (R2) through $\Sigma^*$ at the time $t = \tau$. 

To complete the proof of {\it (b)}, notice that from the construction
of the curves $\Sigma_s$ and $\Sigma_0$ and due to the fact that
$z(\tau) = r'$, $\la(\tau) = r$, we have $z(t) > 0$ and $\la(t) < r$
for $t \in [0,\tau)$.
\medskip

\noindent
Ad {\it (c)}: Assume $\la(0) = r$. Then we have $z_0 \ge r'$.
If $z_0 = r'$, then $\la'(0) = 0$ and $\la''(0) = -z'(0) = \psi_*(x_0) < 0$.
Then, there exists $\tau > 0$ with $\la(\tau) = r$, $\la(t) < r$, $t \in
(0,\tau)$, $\la'(\tau) = 0$, $z(\tau) = r'$. Therefore exists $\sigma \in
(0,\tau)$ with $z(\sigma) > r'$ and $z'(\sigma) = 0$. From the differential
equation for $z$ follows $x(\sigma) > x^*$. Now, note that $(x_0,K_0) \in$
(R3) implies $(x_0,K_{0,+}) \in$ (R3) $\cup$ (R4); further, the solution
does not jump in $(0,\sigma]$ and $u(t) = 1$ in $[0,\sigma]$. Therefore,
from the construction of $\Gamma_2$, we obtain the existence of $\rho \in
(0,\sigma)$ with $(x(\rho),K(\rho)) \in \Sigma^*$ (and even $K(\rho) <
\widetilde{K}_1$). However, since $\la(\rho) < r$, this contradicts Lemma~%
\ref{lem:erreich}. If $z_0 > r'$, then $\la'(0) < 0$ and we obtain a
contradiction analogous to the case $z_0 = r'$.

\noindent
Therefore, $\la_0 < r$. Assume $z_0 \ge 0$. Then $z'(t) = z g(x) - \psi(x)
> 0$, for all $t > 0$ such that $x(t) < \tilde{x}$. Consequently, $z(t) > 0$
as long as $x(t) < \tilde{x}$. We already know that $\la(t) < r$ as long as
$(x(t),K(t)) \in$ (R3) (i.e. no jumps in (R3)). Therefore, we conclude from
the definition of $\Gamma_2$ that $(x,K)$ leaves (R3) through $\Sigma_0$,
i.e. exists $\sigma > 0$ with $(x(\sigma),K(\sigma)) \in \Sigma_0$. From
the definition of $\Sigma_0$, follows $x(\sigma) < \tilde{x}$. Thus, 
$z(\sigma) > 0$ must hold. Now, from $\la_0 < r$, follows the existence of
$\tau > 0$ with $\la(\tau) = r$ and $z(\tau) = r' > 0$, $z'(\tau) =
-\la''(\tau) \ge 0$ (obviously $\tau \ge \sigma$). Then, since $\la(t) = r$
is not allowed in (R3) $\cup\, \Sigma_0\, \cup$ (R2), the case $x(\tau)
< x^*$ can be excluded. Therefore, $x(\tau) = x^*$ must hold, from what
follows $(x(\tau),K(\tau)) \in \Sigma^*$. However, this cannot occur, since
$z(\sigma) > 0$ is not in agreement with Corollary~\ref{kor:nullkurve}.

\noindent
Therefore, $z_0 < 0$. Note that the optimal trajectory meets $\Sigma_0$.
Indeed, this follows from the definition of $\Gamma_4$ and the fact that
$\la(t) < r$ and $z(t) < 0$ as long as $(x(t),K(t)) \in$ (R3).
\medskip

\noindent
Ad {\it (d)}: The case $\la_0 = r$ is excluded arguing as in {\it (c)}.
Assume $z_0 \ge 0$. Then $z'(0) > 0$ and $z'(t) = z g(x) - \psi(x) > 0$,
for all $t > 0$ such that $x(t) < \tilde{x}$. From the differential equation
for $(x,K)$, follows that the solution reaches (R3) with $z(t) > 0$. From
item {\it (c)} we know that this is not possible. Therefore, we have $\la_0
< r$ and $z_0 < 0$. Further, we have $\la(t) < r$ and $z(t) < 0$ as long as
$(x(t),K(t)) \in$ (R4), since otherwise we could repeat the arguments above.
From the construction of $\Gamma_4$, we conclude that the solution meets
$\widetilde{\Sigma}$.
\medskip

\noindent
Ad {\it (e)}: Assume $\la(0) = r$. Then we have $z_0 \ge r'$. \\
If $z_0 = r'$, then $\la'(0) = 0$ and we consider three cases: {\it i)}
$x_0 > x^*$: we have $\la''(0) = \psi_*(x_0) > 0$, contradicting (R);
{\it ii)} $x_0 = x^*$: cannot occur, due to Lemma~\ref{lem:lrplus};
{\it iii)} $x_0 < x^*$: we have $z'(0) > 0$, $\la''(0) < 0$. Since
$\la(t) < r$ for all $t > 0$ is not possible, there exists $\tau > 0$
with $\la(\tau) = r$, $\la(t) < r$, $t \in (0,\tau)$, $\la'(\tau) = 0$,
$z(\tau) = r'$. \\
Then there exists $\sigma \in (0,\tau)$ with $z'(\sigma) = 0$, $z(t) > r'>0$,
$t \in (0,\sigma]$. Consequently, $\psi_*(x(\sigma)) > 0$ and
$x(\sigma) > x^*$. Since $(x_0,K_0) \in$ (R5), then $(x_0,K_{0,+}) \in$ (R5)
too. This fact together with $u(t) = 1$ in $[0,\sigma]$ and the differential
equation for $(x,K)$, implies $K(\sigma) < F(x(\sigma)) / x(\sigma)$,
contradicting Lemma~\ref{lem:sprung1}. If $z_0 > r'$, then $\la'(0) < 0$
and we have again a contradiction.

\noindent
Therefore, $\la_0 < r$. Consequently, there exists $\tau > 0$ with $\la(\tau)
= r$ and $z(\tau) = r' > 0$. Next we exclude two cases: {\it i)} $z_0 < 0$:
there exists $\sigma \in (0,\tau)$ such that $z(\sigma) = 0$, $z'(\sigma)
\ge 0$ and $z(t) < 0$, $t \in [0,\sigma)$. Then $x(\sigma) \le \tilde{x}$
must hold. However, since $u(t) = 0$ in $[0,\sigma)$, we have $x'(t) =
F(x) > 0$, $t \in [0,\sigma)$. Thus we obtain $x(\sigma) > x_0 > \tilde{x}$,
which is a contradiction. {\it ii)} $z_0 = 0$: if $z'(0) < 0$ we obtain a
contradiction arguing as in {\it i)}. If $z'(0) \ge 0$ we have $\psi(x_0)
\le 0$ and $x_0 \le \tilde{x}$, contradicting $(x_0,K_0) \in$ (R5).

\noindent
Therefore, $z_0 > 0$. Finally, we prove that the optimal trajectory meets
$\widetilde{\Sigma}$. We already know that $\la_0 < r$. Then there exists
$\tau > 0$ with $\la(\tau) = r$ and $z(\tau) = r'$. We consider two cases:
{\it i)} $z(\sigma) = 0$, for some $\sigma \in (0,\tau)$: without loss of
generality, we can assume $z'(\sigma) \le 0$. If $z'(\sigma) = 0$, then
$x(\sigma) = \tilde{x}$ and the proof is complete. If $z'(\sigma) < 0$,
then $x(\sigma) > \tilde{x}$ must hold. Moreover, there exists $\rho \in
(\sigma,\tau)$ with $z(\rho) = 0$, $z(t) < 0$ in $(\sigma,\rho)$, $z'(\rho)
\ge 0$. From $x(\sigma) > \tilde{x}$ and $x'(t) = F(x) > 0$, for $t \in
(\sigma,\rho)$, follows $x(\rho) > \tilde{x}$. But this contradicts
$x(\rho) \le \tilde{x}$, which follows from $z(\rho) = 0$, $z'(\rho) \ge 0$.
{\it ii)} $z(t) > 0$ in $[0,\tau]$: since $\la_0 = r$ is not allowed in
(R5), the trajectory must leave (R5) at some time $\sigma \le \tau$. From
$\la(t) < r$ in $[0,\sigma)$ (i.e. no jumps in $[0,\sigma)$), $u(t) = 1$
in $[0,\sigma]$ and the differential equation for $(x,K)$, we conclude that
the trajectory must leave (R5) through $\widetilde{\Sigma}$ as conjectured.
\eop
\end{satz}

\begin{satz}\label{th:gamma}                                         
Let $(x,K,u,\mu)$ be optimal with adjoint variables $z$, $\la$. The
following assertions hold:
\begin{itemize}
\item [(a)] If $(x_0,K_0)$ is on the curve $\Gamma_4$, then we may apply 
Theorem~\ref{th:rrr} (d);
\item [(b)] If $(x_0,K_0)$ is on the curve $\Gamma_3$, then we may apply 
Theorem~\ref{th:rrr} (e).
\end{itemize}
\bew \rm \ Follows along the argumentation above.
\eop
\end{satz}

Resuming, we know that each optimal trajectory meets $\Sigma_0$,
$\widetilde{\Sigma}$ or $\Sigma^*$. We now have to discuss the behavior 
for points on these curves.

\begin{satz}\label{th:sigma0}                                        
Let $(x,K,u,\mu)$ be optimal with adjoint variables $z$, $\la$. The following 
assertions hold:
\begin{itemize}
\item [(a)] Let $(x_0,K_0)$ be on $\Sigma^*$. Then $z_0 =r'$, $\la = r$,
$(x_0,K_{0,+}) = (x^*,K^*)$ and $u \equiv 1$, $\mu = \gamma K dt$;
\item [(b)] Let $(x_0,K_0)$ be on $\Sigma_0$. Then $z_0 = 0$, $\la_0 < r$
and there exists $\tau > 0$ with $(x(\tau),K(\tau)) \in \Sigma^*$ and
$u \equiv 1$, $\mu \equiv 0$ in $(0,\tau)$;
\item [(c)] Let $(x_0,K_0)$ be on $\widetilde{\Sigma}$. Then $z_0 = 0$,
$\la_0 < r$ and there exists $\tau_1 \ge 0$, $\tau_2 > \tau_1$ with
$(x(\tau_1),K(\tau_1)) = (\tilde{x},\ttilde{K})$, $(x(\tau_2), K(\tau_2))
\in \Sigma^*$. Moreover
$$ \mu_{|(0,\tau_1)} \equiv 0,\ \
   u(t) = K(t)^{-1}F(\tilde{x})\tilde{x}^{-1},\ t \in (0,\tau_1),\ \
   \mu_{|(\tau_1,\tau_2)} \equiv 0,\ \
   u_{|(\tau_1,\tau_2)} \equiv 1 . $$
\end{itemize}
\bew \rm \
Note that {\it (a)} was already proved in Lemma~\ref{lem:lrminus}. Item
{\it (b)} is proved exactly in the same way as Theorem~\ref{th:rrr}
(b). Now we prove {\it (c)}. \\
Assume $\la(0) = r$. Then $z_0 \ge r'$. If $z_0=r'$, then $\la'(0) = 0$,
$z'(0) = -\psi_*(\tilde{x}) > 0$, $\la''(0) = -z'(0) < 0$. If $z_0 > r'$,
then $\la'(0) < 0$. In each case there exists $\tau > 0$ with $\la(\tau) = r$,
$\la(t) < r$, $t \in (0,\tau)$, and we obtain $\la'(\tau) = 0$, $z(\tau) = r'$
and $\sigma \in (0,\tau)$ with $z'(\sigma) = 0$. Moreover $\la(t) < r$,
$t \in (0,\tau)$, $z(t) > r'$, $t \in (0,\tau)$. Thus we arrive in (R4)
with $z(t) > 0$, $\la(t) < r$. But this is not in agreement with the results
for domain (R4). Therefore, $\la_0 < r$ must hold. \\
If $z_0 < 0$, then the trajectory meets (R5) with $z(t) < 0$, $\la(t) < r$,
which is not possible due to the results for (R5). If $\la_0 < r$ and $z_0
> 0$, the trajectory reaches (R4) with $z(t) > 0$, $\la(t) < r$, but this is
not allowed due the results for (R4). \\
Thus we have $\la_0 < r$, $z_0 = 0$ and, consequently, $z'(0) = 0$.
Repeating the arguments above, we conclude that the trajectory cannot leave
$\widetilde{\Sigma}$ as long as $K(t) \ge \ttilde{K}$ holds. Therefore,
there exists $\tau_1 > 0$ with $x(t) = \tilde{x}$, $u(t) = F(\tilde{x})
\tilde{x}^{-1} K(t)^{-1}$, $t \in [0,\tau_1]$. In $(\tilde{x},\ttilde{K})$
the trajectory has to follow $\Gamma_2$, otherwise it would enter
(R2) strictly above $\Gamma_2$ and below $\Gamma_3$ with $z \le 0$ 
which is not allowed due to Theorem~\ref{th:rrr}~(b).
Thus the trajectory  meets $\Sigma^*$ at a time $\tau_2 > \tau_1 > 0$.
\eop  
\end{satz}
\bigskip

\noindent
{\large\bf Acknowledgments:} \\[1ex]
The authors would like to thank an anonymous referee for the valuable
comments which helped to improve significantly the quality of the text. \\
A.L. is on leave from Department of Mathematics, Federal University of
St.\,Catarina, 88010-970 Florianopolis, Brazil; his work is supported by
the Austrian Academy of Sciences.
%
%
%
\section{Appendix: maximum principle}

The optimal control problem we want to consider is the following one:
 
\begin{quote}
Minimize $\D J(x_0,K_0;\mu,u) := \int_0^{\infty} e^{-\delta t}r \mu(dt) +
         \int_0^{\infty} e^{-\delta t} \{c - px(t)\} u(t)K(t) dt$ \\
subject to $(u,\mu) \in U_{ad} \times C^*$ and
\begin{eqnarray}
x' & = & F(x) - u(t) K x,\ x(0) = x_0 , \\[1ex] \label{eq:pop1}
dK & = & -\gamma K dt + \mu(dt),\ K(0) = K_0 , \label{eq:kap1}
\end{eqnarray}
\end{quote}
where
\begin{eqnarray*}
U_{ad} & := & \{ v \in L_{\infty}[0,\infty)\ |\ 0 \le v(t) \le 1 \mbox{
              \it a.e. in } [0,\infty) \} , \\
C^*    & := & \{ \mu\ |\ \mu \text{ non--negative Borel measure on }
               [0,\infty) \} .
\end{eqnarray*}
This problem is denoted by $P(x_0,K_0)$. Let $(x,K,u,\mu)$ be a solution of 
the problem. The idea for proving a maximum principle comes from \cite{Ha} 
and \cite{We}. Halkin, however, uses a different solution concept, avoiding 
the use of Bellman's principle to analyze problems with infinite horizon.

Define the so--called Hamilton function $\widetilde{H}$ by
$$ \widetilde{H}(t,\tilde{x},\widetilde{K},w,\tilde{\la}_1,\tilde{\la}_2,\eta)
   := \tilde{\la}_1 (F(\tilde{x}) - w \widetilde{K} \tilde{x})
      - \tilde{\la}_2 \gamma \widetilde{K}
      - \eta e^{-\delta t} (c-p\tilde{x}) w \widetilde{K} . $$
Let $(T_k)_{k\in\N}$ be a sequence in $(0,\infty)$ with $\lim_k T_k = \infty$ 
and consider for each $k \in \N$ the following problem $P_k(x_0, K_0)$:
\begin{quote}
Minimize $\D J(x_0,K_0;\nu,w) := \int_0^{T_k} e^{-\delta t}r \nu(dt) +
         \int_0^{T_k} e^{-\delta t} \{c - py(t)\} w(t) l(t) dt$ \\
subject to $(w,\nu) \in U_{{ad},k} \times C^*_k$ and
\begin{eqnarray}
y' &=& F(y) - w(t) l y,\ y(0) = x_0,\ y(T_k) = x(T_k), \\[1ex] \label{eq:popk}
dl &=& -\gamma l dt + \nu(dt),\ l(0) = K_0,\ l(T_k) = K(T_k), \label{eq:kapk}
\end{eqnarray}
\end{quote}
where
\begin{eqnarray*}
U_{{ad},k} & := & \{ v \in L_{\infty}[0,T_k]\ |\ 0 \le v(t) \le 1 \mbox{
                  \it a.e. in } [0,T_k] \} ,\\
C^*_k      & := & \{ \mu\ |\ \mu \text{ nonnegative Borel measure on }
                  [0, T_k] \} .
\end{eqnarray*}
{\bf Assumption:} $T_k$ is a point of continuity for each $k \in \N$. (Clearly 
such a choice of $(T_k)_{k \in \N} $ is always possible since $K$ possesses 
only a countable number of jumps.)

Once one has proved the Bellman's optimality principle for control problems 
with infinite horizon (see \cite{BL}), one concludes that for each $k \in \N$, 
$(x_{|[0,T_k]}, K_{|[0,T_k]}, \mu_{|[0,T_k]}, u_{|[0,T_k]})$ is a solution of 
$P_k(x_0,K_0)$. With the maximum principle proved in \cite{VP}, we obtain 
$\tilde{\la}_{1,0,k}$, $\tilde{\la}_{2,0,k}$, $\eta_k \in \R$ such that there 
exists $\tilde{\la}_{1,k}$, $\tilde{\la}_{2,k}$ with
\begin{eqnarray*}
\tilde{\la}_{1,0,k}^2 + \tilde{\la}_{2,0,k}^2 + \eta^2_k &=& 1,\
\eta_k \ge 0, \\
x' & = & F(x) - u(t) K x,\ x(0) = x_0 , \\
dK & = & -\gamma K dt + \mu(dt),\ K(0) = K_0, \\
\tilde{\la}_{1,k}' & = & -\tilde{\la}_{1,k}(F'(x) - u(t) K)
                         - \eta_k e^{-\delta t} p u(t) K, \\
\tilde{\la}_{2,k}' & = & \tilde{\la}_{1,k} x u + \gamma \tilde{\la}_{2,k}
                         + \eta_k e^{-\delta t} (c - p x) u(t), \\
\tilde{\la}_{2,k}(t) - \eta_k e^{-\delta t} r & \le & 0
\text{ for all } t \in [0,T_k], \\
\tilde{\la}_{1,k}(0) & = & \tilde{\la}_{1,0,k},\
                           \tilde{\la}_{2,k}(0) = \widetilde{\la}_{2,0,k}, \\
\tilde{\la}_{2,k}(t) - \eta_k e^{-\delta t}r & = & 0 \
                           \mu\text{--a.e. in }[0,T_k], \\
\widetilde{H}(t,x(t),K(t),u(t),\tilde{\la}_{1,k}(t),\tilde{\la}_{2,k}(t), 
              \eta_k) & = & \max_{w \in [0,1]}
\widetilde{H}(t,x(t),K(t),w,\tilde{\la}_{1,k}(t),\tilde{\la}_{2,k}(t),\eta_k)\\
                      & & \hskip6cm \mbox{\it a.e. in } [0,T_k] .
\end{eqnarray*}

Without loss of generality we may assume that the sequences 
$(\tilde{\la}_{1,0,k},\tilde{\la}_{2,0,k})_{k\in\N}$ and $(\eta_k)_{k\in\N}$ 
converge. Let
$$ \tilde{\la}_{1,0} := \lim_{k\to\infty}
   \tilde{\la}_{1,0,k},\ \tilde{\la}_{2,0} := \lim_{k\to\infty}
   \tilde{\la}_{2,0,k},\ \eta := \lim_{k\to\infty} \eta_k . $$
Then, due to the continuous dependence of the solution on initial data 
and parameter (see \cite{Sch}) we obtain
$$ \tilde{\la}_1 := \lim_{k\to\infty} \tilde{\la}_{1,k},\ \tilde{\la}_2 :=
   \lim_{k\to\infty} \tilde{\la}_{2,k} \,, $$
uniformly on each interval $[0,T],T > 0$. This gives the desired maximum
principle for $P(x_0,K_0)$:
\medskip

There exist $\tilde{\la}_{1,0}$, $\tilde{\la}_{2,0}$, $\eta \in \R$
such that there exists $\tilde{\la}_{1}$, $\tilde{\la}_{2}$ with
\begin{eqnarray*}
\tilde{\la}_{1,0}^2 + \tilde{\la}_{2,0}^2 + \eta^2 &\neq & 0,\ \eta \ge 0, \\
x' & = & F(x) - u(t) K x,\ x(0) = x_0, \\
dK & = & -\gamma K dt + \mu(dt),\ K(0) = K_0, \\
\tilde{\la}_1' & = & -\tilde{\la}_1 (F'(x) - u(t)K) - \eta e^{-\delta t}
                     p u(t) K, \\
\tilde{\la}_2' & = & \tilde{\la}_1 x u + \gamma \tilde{\la}_2
                     + \eta e^{-\delta t} (c - p x) u(t), \\
\tilde{\la}_1(0) & = & \tilde{\la}_{1,0},\
\tilde{\la}_2(0) = \tilde{\la}_{2,0}, \\
\tilde{\la}_2(t) - \eta e^{-\delta t} r & \le & 0 \text{ for all }
                   t \in [0,\infty), \\
\tilde{\la}_2(t) - \eta e^{-\delta t} r & = & 0 \
                   \mu \text{--a.e. in } [0,\infty), \\
\widetilde{H}(t,x(t),K(t),u(t),\tilde{\la}_1(t),\tilde{\la}_2(t),\eta) & = &
\max_{w\in[0,1]} \widetilde{H}(t,x(t),K(t),w,\tilde{\la}_1(t),
                          \tilde{\la}_2(t),\eta) \text{ a.e. in } [0,\infty) .
\end{eqnarray*}

\renewcommand{\baselinestretch}{1}
\small \normalsize

\newpage
\pagestyle{plain}
\listoffigures

\newpage

\unitlength1mm

\unitlength1mm

\unitlength1mm

\unitlength1mm

\end{document}